\documentclass[fontsize=12pt,a4paper,headings=normal,
twoside=false,leqno,parskip=half-,abstract=true]{scrartcl}
\usepackage[english]{babel}
\usepackage[utf8]{inputenc}
\usepackage{hyperref}
\hypersetup{
 pdftitle={Periodic Solutions in Symmetric Monotone Feedback},
 pdfauthor={Alejandro López-Nieto},
 colorlinks=true,
 linkcolor=blue,
 citecolor=blue,
 filecolor=blue,
 urlcolor=blue}
\usepackage{xpatch}
\xpatchcmd{\proof}{\itshape}{\normalfont\proofnamefont}{}{}

\newcommand{\proofnamefont}{}

\usepackage{caption}
\usepackage{amsmath,amsthm}
\usepackage{amssymb} 
\usepackage[normalem]{ulem}

\newtheorem{theorem}{Theorem}[section]
\newtheorem{lemma}[theorem]{Lemma}
\newtheorem{corollary}[theorem]{Corollary}

\theoremstyle{definition}

\theoremstyle{remark}
\newtheorem{remark}[theorem]{\emph{Remark}}

\numberwithin{equation}{section}

\usepackage{cellspace}
\usepackage{diagbox}
\usepackage{pgfplots}
\usepackage{float}
\usepackage{makecell}
\usepackage{enumerate}

\newcommand\cincludegraphics[2][]{\raisebox{-0.3\height}{\includegraphics[#1]{#2}}}
\newcommand{\email}[1]{\href{mailto:#1}{#1}}

\definecolor{refkey}{rgb}{1,0,0}
\definecolor{labelkey}{rgb}{1,0,0}
\begin{document}

\title{Global Bifurcation of Periodic Solutions in Delay Equations with Symmetric Monotone Feedback}

\author{Alejandro López-Nieto\thanks{Department of Mathematics, National Taiwan University, No. 1, Sec. 4, Roosevelt Road, 10617 Taipei, Taiwan; \email{alopez@ntu.edu.tw}.}}
\date{ }

\maketitle
\thispagestyle{empty}
\newpage

\begin{abstract}
We study the periodic solutions of the delay differential $\dot{x}(t)=f(x(t),x(t-1))$, where $f$ scalar is strictly monotone in the delayed component and has even-odd symmetry. We completely describe the global bifurcation structure of periodic solutions via a period map originating from planar ordinary differential equations. Moreover, we prove that the first derivative of the period map determines the local stability of the periodic orbits. This article builds on the pioneering work of Kaplan and Yorke, who found some symmetric periodic solutions for $f$ with even-odd symmetry. We enhance their results by proving that all periodic solutions are symmetric if $f$ is in addition monotone.
\end{abstract}

\section{Introduction}
The main goal of this work is to give a complete characterization of the global branches of periodic solutions and their local stability in delay differential equations possessing monotone feedback and even-odd symmetries.

\subsection{Motivation}

We consider the delay differential equation (abbr. DDE)
\begin{align}\label{eq:1.2}
\dot{x}(t)= f(x(t),x(t-1)),
\end{align}
where $f\in C^2(\mathbb{R}^2,\mathbb{R})$ satisfies the \textit{monotone delayed feedback assumption} $\partial_2 f \neq 0$. The DDE \eqref{eq:1.2} has been used widely in the literature to capture time lags arising from the physical constraints of real-life models. For instance, delays appear in population dynamics due to maturation processes that change the role of the species involved. Examples include the Hutchinson equation for rodents \cite{Hu48}, the Lasota--Wazewska equation modelling red blood cell production \cite{LaWa76}, and the Mackey--Glass equation describing circulating blood cells \cite{MaGla77}. In atmospheric sciences, heat transfer lags mediated by oceanic currents reproduce the ENSO phenomenon; see \cite{SuaSch88}. In economics, market reaction lags induce price fluctuations; see \cite{BrErWa04}. In the mammalian segmentation clock \cite{YoKo20}, \eqref{eq:1.2} models oscillations in gene expression due to delays in RNA transcription. Finally, \eqref{eq:1.2} also serves as a technical instrument in describing the motion of colloidal particles governed by a Langevin equation with a memory kernel \cite{LoHerKl21}.

We consider two feedback types. If $\partial_2 f > 0$, then we say that the DDE \eqref{eq:1.2} has \textit{positive monotone feedback}. Otherwise, if $\partial_2f < 0$, then \eqref{eq:1.2} has \textit{negative monotone feedback}. In analogy to Lotka--Volterra systems, the positive and negative monotone feedback types characterize if the interaction between the dynamical variables $x(t)$ and $x(t - 1)$ is cooperative or competitive, respectively. For illustration, we regard $x(t)$ and $x(t - 1)$ as two separate species (the current population and the old one). Then, positive monotone feedback represents a relation ruled by intergenerational cooperation where $x(t)$ benefits from a thriving population $x(t - 1)$ in the past. Conversely, if the feedback is negative, then there is intergenerational competition, that is, the population growth is inhibited by prosperity at a previous time. This is the Malthusian scenario, where a large population triggers an inhibitory response.

We aim to find all \textit{periodic solutions} of the DDE \eqref{eq:1.2}, that is, all nonconstant solutions $x^\ast(t)$ for which there exists a number $p>0$ such that $x^\ast(t+p)=x^\ast(t)$ for all $t\in\mathbb{R}$. Any such $p>0$ is called a \textit{period} of  $x^\ast(t)$ and the smallest period is called the \textit{minimal period} of $x^\ast(t)$. Our reasons to look for periodic solutions are twofold. First, the existence results in \cite{MP88, Nuss73} show that periodic solutions are abundant in \eqref{eq:1.2} and often appear via topological fixed point theorems. Hence, we want to further understand the mathematical structure of the set of periodic solutions. Second, by a Poincar\'{e}--Bendixson theorem \cite{MPSe962}, periodic solutions are one of the cornerstones of the global dynamics of \eqref{eq:1.2}. Therefore, the periodic solutions of \eqref{eq:1.2} act both as stable observable states and as separatrices between regions of qualitatively different dynamical natures. In conclusion, an accurate description of the periodic solutions of \eqref{eq:1.2} is necessary to explain its global dynamics precisely; see \cite{LN23}.

\subsection{Dimension reduction}
A technical complication in the analysis of the DDE \eqref{eq:1.2} is that it generates an infinite-dimensional dynamical system. We denote by $C$ the Banach space $C^0([-1,0],\mathbb{R})$ equipped with the supremum norm. A continuous curve $x:[-1,t^\ast)\to\mathbb{R}$ with $t^\ast>0$ is a \textit{solution} of \eqref{eq:1.2} if $x(t)$ is differentiable and satisfies the differential equality \eqref{eq:1.2} for $t\in [0,t^\ast)$. Here  $\dot{x}(0)$ refers to the right-side derivative at $t=0$. The function $\psi\in C$ such that $x(\theta)=\psi(\theta)$ for $\theta\in[-1,0]$ is called the \textit{initial condition} of the solution. Following \cite{HaLu93}, we regard solutions of \eqref{eq:1.2} as curves of functions $x_t(\theta):= x(t+\theta)$ for $\theta\in[-1,0]$. The family of transformations $S(t):C\to C$ defined via $S(t)x_0:= x_t$ is called the \textit{semiflow} of \eqref{eq:1.2}.

We highlight three widespread ways to prove the existence of periodic solutions in the DDE \eqref{eq:1.2}. The first approach is the Hopf bifurcation of periodic orbits from a known reference equilibrium. This method only requires us to study changes in the spectrum of a linear operator as a parameter varies. In the case of \eqref{eq:1.2}, the problem has been well studied \cite{HaLu93} and it reduces to finding purely imaginary solutions of the transcendental characteristic equation
\begin{align}
\mu=A+B\mathrm{e}^{-\mu}\quad\text{for } A,B\in\mathbb{R},\quad B\neq 0.
\end{align}
However, the bifurcation is only local, meaning that the solutions obtained in this way stay confined near the reference equilibrium.

The second approach relies on homotopy-invariant indexes to achieve periodic solutions by a continuation process from better understood equations; see \cite{MP88, Nuss73}. The main drawback is that any two different equations within the same homotopy class are seen as equal by this method. The outcome is a loss of information on the number of periodic solutions and their amplitudes.
 
The third approach is to use symmetry. Kaplan and Yorke \cite{KaYo74} considered the delay-only DDE
\begin{align}\label{eq:1.3}
\dot{x}(t)=-g(x(t-1)).
\end{align}
Here, $g$ is odd and has the opposite sign of $x(t-1)$, that is, 
\begin{align}
g(-{u})=-g({u})\quad \text{and satisfies}\quad {u} g({u})>0\;\text{for}\;{u}\neq 0.
\end{align}
If we further assume the \textit{symmetry Ansatz}
\begin{align}\label{eq:1.4}
x^\ast(t-2)=-x^\ast(t),
\end{align} 
then, setting $({u_1}(t),{u_2}(t)):=(x^\ast(t),x^\ast(t-1))$, we obtain that $({u}_1(t), u_2(t))$ solves the ODEs in $\mathbb{R}^2$
\begin{align}\label{eq:1.5}
\begin{split}
\dot{u}_1&=-g({u_2}),\\
\dot{u}_2&=g({u_1}).
\end{split}
\end{align}
Hence, any periodic solution of the DDE \eqref{eq:1.2} with the symmetry \eqref{eq:1.4} yields a solution of the ODEs \eqref{eq:1.5}. Conversely, the ODEs \eqref{eq:1.5} are Hamiltonian and preserve the energy function
\begin{align}
E({u_1},{u_2}):=\int_0^{u_1} g(s)\,\mathrm{d}s+\int_0^{u_2} g(s)\,\mathrm{d}s.
\end{align}
The level sets of $E$ are the trajectories of periodic solutions $({u}_1^\ast(t),{u}_2^\ast(t))$ winding around the only equilibrium of \eqref{eq:1.5} at $(0,0)$ and have minimal period $p>0$. Moreover, denoting $G({u_1},{u_2}):=(-g({u_2}),g({u_1}))$, the ODEs \eqref{eq:1.5} are \textit{equivariant} with respect to the rotation by $\pi/2$ given by
\begin{align}\label{eq:1.7}
\varrho({u_1},{u_2}):=({u_2},-{u_1}),\quad \text{that is},\quad\varrho^{-1}\circ G\circ \varrho=G.
\end{align}
Additionally, \eqref{eq:1.5} are \textit{reversible} with respect to the reflections by $\sigma$ and $-\sigma$ defined by
\begin{align}\label{eq:1.8}
\pm\sigma({u_1},{u_2}):=\pm({u_2},{u_1}),\qquad \text{that is},\qquad\left(\pm\sigma^{-1}\right)\circ G\circ \left(\pm\sigma\right)=-G.
\end{align}

Kaplan and Yorke \cite{KaYo74} showed that the equivariance \eqref{eq:1.7} implies the delay relation
\begin{align}
    {u}_2^\ast(t)={u}_1^\ast\left(t-\frac{(4n-3)p}{4}\right),\quad \text{for some } n\in \mathbb{N}.
\end{align}
In this way, the first component $u_1(t)$ of any periodic solution of \eqref{eq:1.5} with minimal period $p=4/(4n-3)$ yields a periodic solution of the DDE \eqref{eq:1.2}. Furthermore, any periodic solution of the DDE \eqref{eq:1.2} obtained in this way automatically satisfies the symmetry Ansatz \eqref{eq:1.4}.

We highlight that the monotone feedback assumption $g'({u})\neq 0$ is not needed for the symmetry approach to work. Thus, the Hamiltonian ODEs \eqref{eq:1.8} only find periodic solutions of the DDE \eqref{eq:1.3} that satisfy the symmetry Ansatz \eqref{eq:1.4}. That is, the symmetry Ansatz \eqref{eq:1.4} is useful for studying the existence and dynamical properties of periodic solutions, \textit{but it does not capture all the possible periodic solutions}. A counterexample is the solutions breaking the symmetry \eqref{eq:1.4} that appear by the period doubling bifurcations in \cite{DorLa95}. Hence, in the absence of monotone feedback,
there are high-dimensional dynamics involved in the formation process of periodic solutions in the DDE \eqref{eq:1.3}.

In this article, we consider a similar setting to that of Kaplan and Yorke. However, in addition, we assume that $g$ is a monotone function. Under this further assumption, all periodic solutions of the DDE \eqref{eq:1.3} satisfy the symmetry Ansatz \eqref{eq:1.4} and therefore yield solutions of the two-dimensional ODEs \eqref{eq:1.5}. The main tool for bridging the gap between the infinite-dimensional DDE \eqref{eq:1.3} and the two-dimensional ODEs \eqref{eq:1.5} is an $\mathbb{R}^2$-embedding theorem for periodic orbits of DDEs with monotone feedback \cite{MPSe962}. Furthermore, we refine the spectral analysis in \cite{MPNu13} for a precise quantification of the local stability of the periodic solutions.

\subsection{Main results}
We say that a function $f$ has \textit{even-odd symmetry} if
\begin{equation}\label{even-odd-symmetry}
    f({u_1},{u_2})=f(-{u_1},{u_2}) \quad \text{and}\quad f({u_1},{u_2})=-f({u_1},-{u_2}).
\end{equation}
This setting is more general than Kaplan and Yorke's since it allows instantaneous feedback, that is, $g$ in the DDE \eqref{eq:1.3} can depend on $x(t)$. Denoting by $\partial_j$ the partial derivative in the $j$-th component, we define the set of nonlinearities with \textit{symmetric positive monotone feedback} by
\begin{align}\label{eq:1.11}
\mathfrak{X}^+:= \left\{f\in C^2(\mathbb{R}^2,\mathbb{R}) : \begin{array}{l}
f \text{ has the symmetry \eqref{even-odd-symmetry} and}\\\partial_2  f({u_1},{u_2})>0\text{ for all } ({u_1},{u_2})\in\mathbb{R}^2
\end{array}\right\}.
\end{align}
Analogously, the set of nonlinearities with \textit{symmetric negative monotone feedback} is $\mathfrak{X}^-:=-\mathfrak{X}^+$. Finally, the set of nonlinearities with \textit{symmetric monotone feedback} is the disjoint union $\mathfrak{X}:=\mathfrak{X}^+\cup\mathfrak{X}^-.$

\begin{remark}
    We say that the functions in $\mathfrak{X}$ have symmetric monotone feedback because, analogously to the Kaplan--Yorke scenario \eqref{eq:1.3}, the feedback strength $|f(x(t), x(t - 1))|$ grows with the signal strength given by $|x(t)|$ and $|x(t - 1)|$. In particular, the feedback strength is invariant under sign changes of the signal.
\end{remark}

In this article, we give a full description of the set of periodic solutions of the DDE \eqref{eq:1.2} with $f\in \mathfrak{X}$. We address two aspects: the global existence and local stability of periodic solutions. There are three main advantages to our approach:
\begin{enumerate}
\item We find all periodic solutions of the DDE \eqref{eq:1.2}, as long as a yet to be defined period map is not constant on any interval.
\item We quantify precisely the local stability of periodic solutions.
\item Our approach is constructive and visual, and it explains the formation mechanism of periodic solutions.
\end{enumerate}

For the existence part, we extend the idea of Kaplan and Yorke in \cite{KaYo74}. We consider the more general DDE \eqref{eq:1.2} with $f\in\mathfrak{X}$ and study periodic solutions $x^\ast(t)$ satisfying the symmetry Ansatz \eqref{eq:1.4}. By even-odd symmetry, the curve $({u_1}(t),{u_2}(t)):=(x^\ast(t),x^\ast(t-1))$ solves the ODEs
\begin{align}\label{eq:1.12}
\begin{split}
\dot{u}_1&=f({u_1},{u_2}),\\
\dot{u}_2&=-f({u_2},{u_1}).
\end{split}
\end{align}
Here $f\in\mathfrak{X}$ ensures that the vector field $F({u_1},{u_2}):=(f({u_1},{u_2}),-f({u_1},{u_2}))$ is equivariant under the rotation \eqref{eq:1.7} and is reversible under the reflections \eqref{eq:1.8}.

These symmetry properties guarantee that the integral curves of the planar ODEs \eqref{eq:1.12} are closed curves centered around the origin $(0, 0)$ and foliate all $\mathbb{R}^2$. Unlike in \cite{KaYo74}, \eqref{eq:1.12} are not necessarily the Hamiltonian vector field of an energy function $E$. We describe the solutions of \eqref{eq:1.12} via the following lemma.

\begin{lemma}[period map]\label{Lemma1.1}
Consider the ODEs \eqref{eq:1.12} with $f\in\mathfrak{X}$.  Then the only equilibrium of \eqref{eq:1.12} is $(0,0)$. Moreover, the solution $({u_1}(t; a),{u_2}(t; a))$ with initial condition $({u_1}(0; a),{u_2}(0; a))=(a, 0)$, for $a>0$, is periodic with minimal period ${p}_f(a)$.
Furthermore, if $f\in C^k(\mathbb{R}^2,\mathbb{R})$ (resp., $f$ is analytic), $k \geq 2$, then the map ${p}_f:(0,\infty)\to \mathbb{R}$ taking the amplitude $a>0$ to the minimal period ${p}_f(a)$ of $({u_1}(t; a),{u_2}(t; a))$ is $C^k$ (resp., analytic).
\end{lemma}

The function ${p}_f$ defined Lemma \ref{Lemma1.1} is called the \textit{period map} of $f\in \mathfrak{X}$. We say that ${p}_f$ is \textit{locally nonconstant at $a>0$} if ${p}_f$ is not constant on any interval $I$ such that $a\in I\subset(0,\infty)$. In particular, ${p}_f$ is called \textit{locally nonconstant} if it is locally nonconstant at all $a>0$.

\begin{figure}
\includegraphics[width=\textwidth]{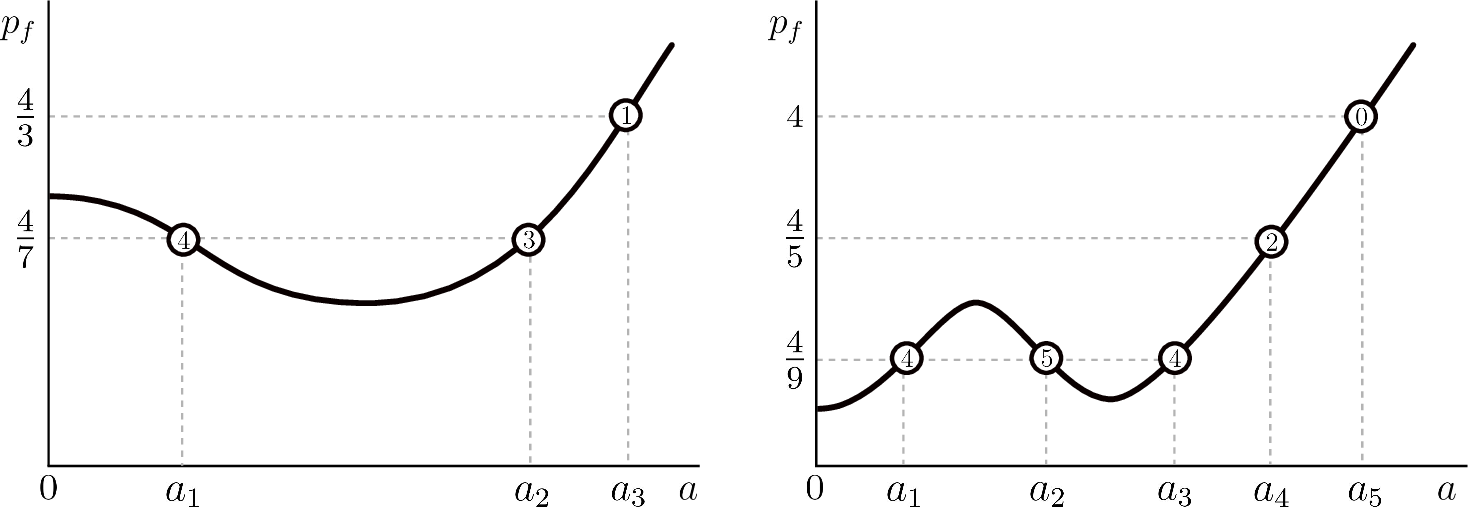}
\caption{(Left) Schematic graph of ${p}_f$ for a nonlinearity $f\in\mathfrak{X^+}$, by Theorem \ref{Theorem1.2}, the DDE \eqref{eq:1.2} possesses solely three periodic solutions with amplitudes $a_i$, marked by circles, the number inside each circle indicates the unstable dimension given by Theorem \ref{Theorem1.3}. (Right) Schematic graph of ${p}_f$ for a nonlinearity $f\in\mathfrak{X^-}$, by Theorem \ref{Theorem1.2}, the DDE \eqref{eq:1.2} possesses a total of five periodic solutions, marked by circles, the number inside each circle indicates the unstable dimension given by Theorem \ref{Theorem1.3}.}
\label{fig:1}
\end{figure}

\begin{theorem}[characterization of periodic solutions]\label{Theorem1.2}
Consider the DDE \eqref{eq:1.2} with nonlinearity $f\in\mathfrak{X}$ such that its period map ${p}_f$ is locally nonconstant. Then any periodic solution $x^\ast(t)$ of \eqref{eq:1.2} has (not necessarily minimal) period 4 and satisfies the symmetry Ansatz \eqref{eq:1.4}. More precisely, let $\bar{x}:= \max_{t\in[0,4]}x^\ast(t)$ be the amplitude of $x^\ast(t)$. Then the minimal period of $x^\ast(t)$ is ${p}_f(\bar{x})$ and it satisfies:
\begin{align}\label{eq:1.14}
{p}_f(\bar{x})=\frac{4}{4n-1}, \quad \text{for some}\;n\in\mathbb{N}\text{, if } f\in\mathfrak{X}^+,
\end{align}
\begin{align}\label{eq:1.15}
{p}_f(\bar{x})=\frac{4}{4n-3}, \quad \text{for some}\;n\in\mathbb{N}\text{, if }f\in\mathfrak{X}^-.
\end{align}

Conversely, let $({u}_1^\ast(t),{u}_2^\ast(t))$ be a periodic solution of the ODEs \eqref{eq:1.12} with amplitude $a:=\max_{t\in\mathbb{R}}{{u}_1^\ast(t)}$. If $f\in\mathfrak{X}^+$ (resp., $f\in\mathfrak{X}^-$) and ${p}_f(a)=4/(4n-1)$ (resp., ${p}_f(a)=4/(4n-3)$), then ${u}_1^\ast(t)$ is a periodic solution of the DDE \eqref{eq:1.2}.
\end{theorem}
\begin{remark}
    Theorem \ref{Theorem1.2} gives sufficient conditions such that all periodic solutions of the DDE \eqref{eq:1.2} satisfy the symmetry Ansatz \eqref{eq:1.4}. Hence, the infinite-dimensional periodic solutions of the DDE appear as solutions of a boundary value problem for planar ODEs. Visually, the amplitudes of the periodic solutions of the DDE \eqref{eq:1.2} are the intersections of the period map with the values given by formulas \eqref{eq:1.14}--\eqref{eq:1.15}; see Figure \ref{fig:1}.
\end{remark}

To quantify the local stability of the periodic solutions obtained in Theorem \ref{Theorem1.2}, we define the \textit{orbit} of a periodic solution $x^\ast(t)$ of the DDE \eqref{eq:1.2} with period $p>0$ to be the set $\gamma(x_0^\ast):= \{x^\ast_t : t\in \mathbb{R}\}\subset C$. We denote $\gamma^\ast:=\gamma(x_0^\ast)$, since $\gamma^\ast$ is invariant under the solution semiflow $S(t)$, its local stability is determined by the spectrum of the linearization of the \textit{time-$p$ map} $S(p)$ at any point of $\gamma^\ast$. The \textit{monodromy operator} ${L}$ of $\gamma^\ast$ is defined as
\begin{align}\label{eq:1.9}
{L}:= D_{x^\ast_0}S(p):C\to C.
\end{align}
Here $D_{x^\ast_0}$ denotes the Fr\'{e}chet derivative at $x^\ast_0$. It is well-known that there exists an $n\in\mathbb{N}$ such that ${L}^n$ is a compact operator on $C$. As a consequence, the spectrum $\mathrm{Spec}({L})$ consists solely of eigenvalues with finite algebraic multiplicity called \textit{Floquet multipliers}, and $0$. 

Following \cite{HaLu93}, ${L}$ coincides with the time-$p$ solution operator, with initial time $t=0$, of the linear nonautonomous equation
\begin{align}\label{eq:1.10}
\dot{y}(t)=\partial_1 f(x^\ast(t),x^\ast(t-1))y(t)+\partial_2 f(x^\ast(t),x^\ast(t-1))y(t-1).
\end{align}
In particular, $\dot{x}^\ast(t)$ is a solution of the linearized DDE \eqref{eq:1.10} with period $p$, which implies that $1\in\mathrm{Spec}({L})$. The periodic orbit $\gamma^\ast$ is called \textit{hyperbolic} if the \textit{trivial Floquet multiplier} $1\in\mathrm{Spec}({L})$ has algebraic multiplicity one and no other Floquet multipliers lie on the complex unit circle.

We characterize the local stability of the periodic orbit $\gamma^\ast$ in terms of its Floquet multipliers. More precisely, $\gamma^\ast$ is called \textit{asymptotically stable} (resp., \textit{unstable}) if it is hyperbolic and all of its Floquet multipliers, except the trivial eigenvalue $1$, lie inside the closed complex unit disk (resp., if any of its Floquet multipliers lies outside the closed complex unit disk).
Our measure of the stability of the periodic orbit $\gamma^\ast$ is the \textit{unstable dimension}, denoted by $i(\gamma^\ast)$. Indeed, we define $i(\gamma^\ast)$ as the number of Floquet multipliers of $\gamma^\ast$ lying outside the closed complex unit disk, counting multiplicities. By definition, $i(\gamma^\ast)>0$ if and only if $\gamma^\ast$ is unstable, conversely, if $\gamma^\ast$ is hyperbolic, then it is asymptotically stable if and only if $i(\gamma^\ast)=0$.

\begin{theorem}[local stability of periodic solutions]\label{Theorem1.3}
Consider a periodic solution $x^\ast(t)$ of the DDE \eqref{eq:1.2} with nonlinearity $f\in\mathfrak{X}$. Let there be $n\in \mathbb{N}$ such that $x^\ast(t)$ has amplitude $\bar{x}:= \max_{t\in[0,4]}x^\ast(t)$ and minimal period
\begin{align}
\label{eq:1.16}
{p}_f(\bar{x})=
\begin{cases}
\dfrac{4}{4n-1},\quad \text{if} \quad f\in \mathfrak{X}^+,\\[2ex]
\dfrac{4}{4n-3},\quad \text{if} \quad f\in \mathfrak{X}^-.\
\end{cases}
\end{align}
Then $\gamma^\ast$ is hyperbolic if and only if ${p}_f'(\bar{x})\neq 0$. Moreover, the unstable dimension $i(\gamma^\ast)$ satisfies:
\begin{align}\label{eq:1.18}
i(\gamma^\ast)&=
\begin{cases}
2n-1, &\text{if}\quad {p}_f'(\bar{x})\geq 0,\\
2n, &\text{if}\quad {p}_f'(\bar{x})<0,
\end{cases}\quad\text{and }f\in\mathfrak{X}^+.\\
\label{eq:1.19}
i(\gamma^\ast)&=
\begin{cases}
2n-2, &\text{if}\quad {p}_f'(\bar{x})\geq 0,\\
2n-1, &\text{if}\quad {p}_f'(\bar{x})<0,
\end{cases}\quad\text{and }f\in\mathfrak{X}^-.
\end{align}
\end{theorem}

\begin{remark}
    Theorem \ref{Theorem1.3} relates the unstable dimension of a periodic solution to the period map. In particular, the periodic solutions with a smaller minimal period, and therefore a higher frequency, have a higher unstable dimension than those with minimal period $4$. Moreover, the only asymptotically stable periodic orbits are those whose amplitude $\bar{x}$ satisfies $p_f(\bar{x}) = 4$ and $p_f'(\bar{x}) > 0$. Our results are based on \cite[Section 5]{MPNu13}, which relates the unstable dimension of a periodic solution to the number of sign changes over a unit length interval. In particular, we have enhanced \cite[Theorem 5.1]{MPNu13} by providing the exact value of the unstable dimension in terms of the maps $p_f$ and $p_f'$; see Figure \ref{fig:1}.
\end{remark}

\subsection{Global branches of periodic orbits} An alternative interpretation of Theorems \ref{Theorem1.2}--\ref{Theorem1.3} is to consider the maximal continuation of periodic solutions in the DDE 
\begin{align}\label{eq:1.15b}
\dot{x}(t)={r} f(x(t),x(t-1)).
\end{align}
Here ${r} \neq 0$ is a parameter that acts as the timescale of the DDE \eqref{eq:1.15b}. Let us fix $f\in \mathfrak{X}^+$ such that ${p}_f$ is locally nonconstant. Then Theorem \ref{Theorem1.2} gives a complete description of the branches of periodic orbits in the $({r}, a)$-plane.

\begin{figure}
\includegraphics[width=\textwidth]{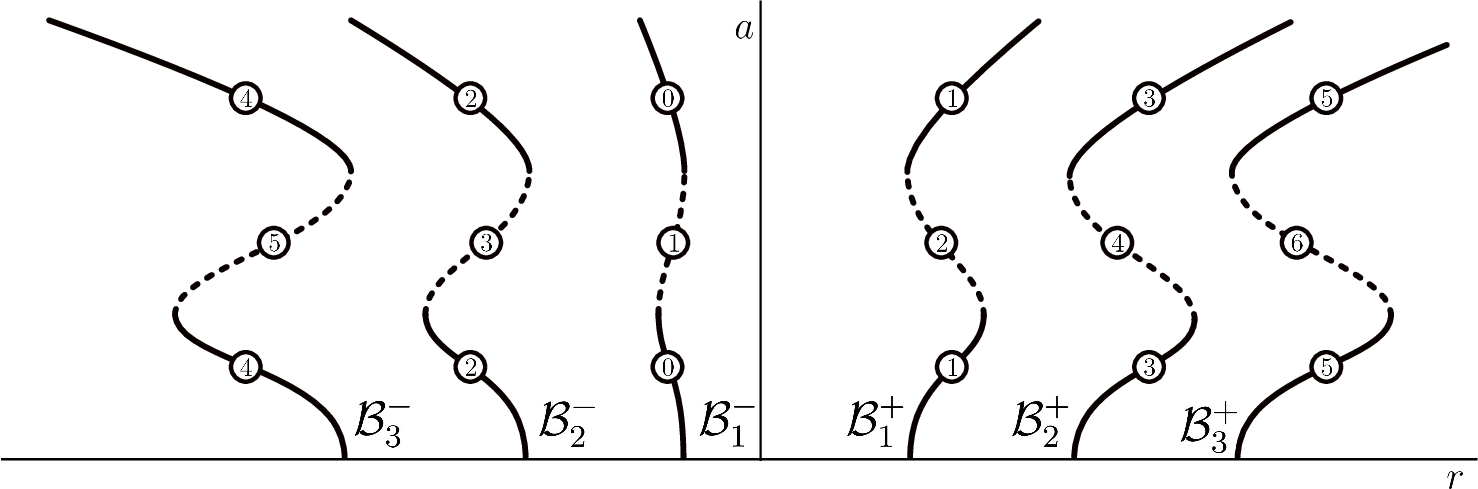}
\caption{First branches of periodic orbits in the DDE \eqref{eq:1.15b} with $f\in\mathfrak{X}^+$. The unstable dimension is indicated in circles. The dashed curves indicate the higher unstable dimension out of the two possible values on each branch. The changes from dashed to solid correspond to saddle-node bifurcations of periodic orbits and the points of intersection with the horizontal axis are Hopf bifurcations.}
\label{fig:0}
\end{figure}

Indeed, let $x^\ast(t)$ be a periodic solution of \eqref{eq:1.15b} for the parameter value ${r}^\ast>0$. We denote the amplitude by $\bar{x}:=\max_{t\in[0,4]} x^\ast(t)$ and the orbit by $\gamma^\ast:= \gamma(x_0^\ast)$. By Theorem \ref{Theorem1.2}, the periodic orbits of \eqref{eq:1.15b} are in one-to-one correspondence with the orbits of the ODEs \eqref{eq:1.12}. More accurately, we identify $\gamma^\ast$ with the point $({r}^\ast, \bar{x})$ in the $({r},a)$-plane.  Furthermore, the parameter ${r}$ in \eqref{eq:1.15b} appears as a time rescaling in the ODEs \eqref{eq:1.12} and produces the identity ${p}_f(\bar{x})/{r}={p}_{{r} f}(\bar{x})$. By the formula \eqref{eq:1.14}, we obtain that
\begin{align}\label{eq:1.15c}
\frac{{p}_f(\bar{x})}{{r}^\ast}={p}_{{r}^\ast f}(\bar{x})=\frac{4}{4n-1},\quad \text{for some } n\in\mathbb{N},
\end{align} 
and $({r}^\ast, \bar{x})$ lies on the \textit{branch} $\mathcal{B}^+_n$ of periodic orbits in the $({r},a)$-plane given by 
\begin{align}
\mathcal{B}^+_n:=\left\{({r},a) \in \mathbb{R}^2 : {r}=\frac{(4n-1){p}_f(a)}{4}\right\} ,\quad \text{for }n\in \mathbb{N}.    
\end{align}
Analogously, we may choose $r^\ast < 0$ above so that $r^\ast f\in\mathfrak{X}^-$. Since the period map is always positive, the formula \eqref{eq:1.15} yields the identity
\begin{align}
-\frac{{p}_f(\bar{x})}{{r}^\ast}={p}_{{r}^\ast f}(\bar{x})=\frac{4}{4n-3},\quad \text{for some } n\in\mathbb{N},
\end{align} 
which produces the branches 
\begin{align}
    \mathcal{B}^-_n:=\left\{({r},a)\in \mathbb{R}^2 : {r}= - \frac{(4n-3){p}_f(a)}{4}\right\},\quad \text{for }n\in \mathbb{N}.
\end{align} An plot of branches for $n = 1, 2,3$ can be seen in Figure \ref{fig:0}.  

The global picture of the bifurcation branches is complemented by two local bifurcation mechanisms that change the local stability of the periodic orbits: Hopf bifurcation and saddle-node bifurcation.

\begin{figure}
\includegraphics[width=\textwidth]{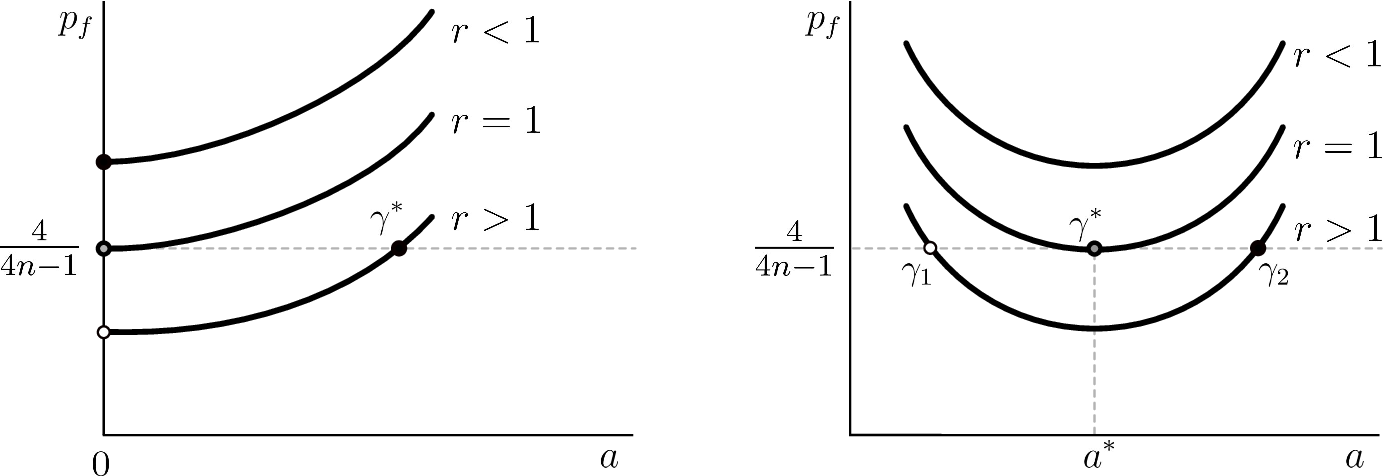}
\caption{(Left) A supercritical Hopf bifurcation in the DDE \eqref{eq:1.15} where $f\in\mathfrak{X}^+$ and  ${p}_f''(0)>0$. At parameter value ${r}=1$, a periodic orbit $\gamma^\ast$ is emitted by the origin and increases its amplitude as ${r}$ grows.  The unstable dimension of $\gamma^\ast$ is $2n-1$. (Right) A saddle-node bifurcation of periodic orbits in the DDE \eqref{eq:1.15} where $f\in\mathfrak{X}^+$ and  ${p}_f''(a^\ast)>0$. At ${r}=1$, the family $\gamma_1$ of  periodic orbits with unstable dimension $2n$ collides with the family $\gamma_2$ with unstable dimension $2n-1$ forming a nonhyperbolic periodic solution $\gamma^\ast$ which has unstable dimension $2n-1$. For ${r}<1$, no periodic solutions exist near the amplitude $a^\ast$ of $\gamma^\ast$.}
\label{fig:4}
\end{figure}

The DDE \eqref{eq:1.15b} undergoes a Hopf bifurcation at ${r}=1$ if ${p}_f(0):=\lim_{a\to 0}{p}_f(a)=4/(4n-1)$ for some $n\in \mathbb{N}$. Depending on the sign of ${p}_f''(0):=\lim_{a\to 0}{p}_f''(a)$,  if ${p}_f''(0)>0$ (resp., ${p}_f''(0)<0$) the bifurcation is \textit{supercritical} (resp., \textit{subcritical}). Thus, a periodic orbit $\gamma^\ast$ with unstable dimension $i(\gamma^\ast)=2n-1$ (resp., $i(\gamma^\ast)=2n$) approaches the equilibrium as ${r}> 1$ decreases (resp., ${r}< 1$ increases). The periodic solution touches $0$ at ${r}=1$, and vanishes for ${r}<1$ (resp., ${r}>1$); see Figure \ref{fig:4} (Left).

Analogously, the DDE \eqref{eq:1.15b} undergoes a saddle-node bifurcation at $r = 1$ if there is an amplitude $a^\ast>0$ for which ${p}_f(a^\ast)=4/(4n-1)$, ${p}_f'(a^\ast)=0$, and ${p}_f''(a^\ast) \neq 0$. At $r = 1$, there exists a single nonhyperbolic periodic orbit $\gamma^\ast$ of \eqref{eq:1.15b} with amplitude $a^\ast$. The orbit $\gamma^\ast$ is the collision point of two families of periodic orbits, $\gamma_1$ and $\gamma_2$, where we set the amplitudes of $\gamma_1$ to be smaller than the amplitudes of $\gamma_2$. If ${p}_f''(a^\ast)>0$ (resp., ${p}_f''(a^\ast)<0$), then $\gamma_1$ and $\gamma_2$ approach one another as ${r}>1$ decreases (resp., ${r}<1$ increases). At ${r}=1$, the families $\gamma_1$ and $\gamma_2$ merge at $\gamma^\ast$ and then vanish for ${r}<1$ (resp. ${r}>1$); see Figure \ref{fig:4} (Right). 

\subsection{Realizability and prototypical example}

We point out that Theorems \ref{Theorem1.2}--\ref{Theorem1.3}, written in terms of ${p}_f$, hide that the period map is, in general, unknown. Since the shape of the period map determines the shape of the branches of periodic orbits in Figure \ref{fig:0}, one crucial question is which functions $p: (0,\infty) \to (0,\infty)$ are realizable as the period map $p_f$ of a nonlinearity $f \in \mathfrak{X}$. We provide an example showing that the shape of $p$ can be prescribed at will, sufficiently close to $a=0$. The so-called \textit{enharmonic oscillator} is the DDE 
\begin{align}\label{eq:1.22b}
    \dot x(t) = 2\pi\Omega\left(x(t)^2 + x(t - 1)^2\right) x(t - 1);
\end{align}
see \cite{LN23}. Here $\Omega: [0, \infty) \to (0,\infty)$ is a positive \textit{frequency function} and the positive monotone feedback condition \eqref{eq:1.11} is equivalent to requiring
\begin{align}\label{eq:1.22c}
    \frac{\Omega'(a^2)}{\Omega(a^2)} >  - \frac{1}{2a^2}, \quad \text{for all }a>0.
\end{align}
Moreover, the ODEs \eqref{eq:1.12} associated to \eqref{eq:1.22b} are
\begin{align}\label{eq:1.22d}
    \begin{split}
        \dot{u}_1 &= 2\pi\Omega\left({u}_1^2 + {u}_2^2\right) {u_2}\\
        \dot{u}_2 &= -2\pi\Omega\left({u}_1^2 + {u}_2^2\right) {u_1}.
    \end{split}
\end{align}
By direct substitution, all solutions of \eqref{eq:1.22d} are of the form
\begin{align}\label{eq:1.22e}
    ({u_1}(t; a), {u_2}(t; a)) = \left(a\sin\left(2\pi\Omega\left(a^2\right)t + t_0\right), a\cos\left(2\pi\Omega\left(a^2\right)t + t_0\right)\right),
\end{align}
with $a \geq 0$ and $t_0 \in \mathbb{R}$. Notice that, first, by \eqref{eq:1.22e}, the period map of \eqref{eq:1.22b} is $p_f(a) = 1/\Omega(a^2)$. Second, the lower bound \eqref{eq:1.22c} goes to minus infinity for sufficiently small amplitudes $0 < a \ll 1$. Hence, the period map of \eqref{eq:1.22b} can have an arbitrary shape sufficiently close to $a = 0$.

\subsection{Connection to the literature}
Our approach based on period maps is inspired by the  bifurcation problem of rotating waves for scalar reaction-diffusion partial differential equations (abbr. PDEs) on the circle. A characterization result for periodic solutions in PDEs that resembles Theorems \ref{Theorem1.2}--\ref{Theorem1.3} was achieved in \cite{FieRoWol04}. However, the period map arises more naturally in PDEs than in the DDE \eqref{eq:1.2} because rotating waves are periodic solutions of planar ODEs.

Further important nonlinearities $f\in\mathfrak{X}^-$ for which the period map is well-understood are \emph{soft-spring nonlinearities}, that is, functions such that $f({u_1},{u_2})/{u_2}=g({u_2})/{u_2}$ is a decreasing function, e.g., 
\begin{align}\label{eq:1.24}
\dot{x}(t)=-{r}\arctan(x(t-1)), \quad {r}>0.
\end{align}
Under these conditions the period map ${p}_f$ is increasing \cite[Theorem 1.3]{Nuss79}. Therefore, Theorem \ref{Theorem1.2} shows that all periodic solutions of the DDE \eqref{eq:1.24} arise by a sequence of Hopf bifurcations from the only equilibrium at $x(t) \equiv 0$. Increasing the parameter ${r}$ in \eqref{eq:1.24} results in new periodic orbits with ever decreasing minimal periods. Our study recovers and improves on previously known results involving the existence, local stability, and uniqueness of stable slowly oscillating periodic solutions \cite{ChWa88,IvLaWa06,KaYo75,Nuss79}. 

The global dynamics of the DDE with positive monotone feedback 
\begin{align}\label{eq:1.25}
\dot{x}(t)=-{r}_1 {r}_2 x(t) +{r}_2\arctan(x(t-1)),\quad {r}_1,{r}_2>0
\end{align} have been widely studied in the literature. In \cite{Kr08, KrWa01, KrWaWu99} the authors give a precise description of the structure of a finite-dimensional geometrical object that attracts all solutions of the DDE \eqref{eq:1.25}. Our results suggest that, in the limit ${r}_1\to 0$, the construction of this global attracting set is explained by a sequence of Hopf bifurcations. A description of an alternative period function associated to \eqref{eq:1.25} via properties of its periodic solutions has been given in \cite{KrGa11}. We point out that these results are complementary to ours in the sense that we use a period map to describe the periodic solutions. The main difference is that in \cite{KrGa11} the period function is ${r}_2$-dependent and therefore collapses near bifurcations points. In turn, we consider an amplitude-dependent function which is defined globally, at the expense of requiring even-odd symmetry of $f$.

The article is structured as follows, in Section \ref{S2} we follow \cite{MPNu13,MPSe962,MPSe961} and introduce the \textit{zero number} and the \textit{Poincar\'{e}--Bendixson property} for DDEs with monotone feedback in Lemmata \ref{Lemma2.1}--\ref{Lemma2.2}. Sections \ref{S3} and \ref{S4} contain the proofs of Lemma \ref{Lemma3.1} and Theorem \ref{Theorem1.2}, and Theorem \ref{Theorem1.3}, respectively.

\section{Zero number and oscillation theory}\label{S2}

In this section, we present the main tools needed to prove our results. First, we introduce the \textit{zero number}, an integer-valued function that measures the frequency of oscillation around zero of a solution to the DDE \eqref{eq:1.2}. For $f\in\mathfrak{X}$, the zero number of a solution is a nonincreasing function in time. This property justifies referring to the zero number as a discrete Lyapunov function and endows the global dynamics of \eqref{eq:1.2} with an explicit {Morse decomposition}; see \cite{MP88}. The zero number has an analogue in scalar PDEs \cite{An88, Ma82}, where it plays a fundamental role to describe the global dynamics \cite{FieRoWol04, FuRo91}. Our use of the zero number for DDEs follows \cite{MPSe962,MPSe961}. At the end of the section, we follow \cite{MPNu13} to discuss a Sturm--Liouville-like oscillation theory induced by the zero number at the linear level. 

Given a function $\psi\in C\setminus\{0\}$, the \textit{sign changes of }$\psi$ are given by
\begin{align}\label{eq:2.1}
\mathrm{sc}(\psi):= \sup\{k\in\mathbb{N} : \psi(\tau_i)\psi(\tau_{i+1})<0,\;-1<\tau_1<\dots<\tau_{k+1}<0\}.
\end{align}
For $\psi\in C\setminus\{0\}$ and $\mathrm{sc}(\psi)<\infty$, we define the \textit{zero number} functions
\begin{align}\label{eq:2.2}
z^+(\psi):= \begin{cases}
\mathrm{sc}(\psi),&\text{ if }\mathrm{sc}(\psi)\text{ is even,}\\
\mathrm{sc}(\psi)+1,&\text{ if }\mathrm{sc}(\psi)\text{ is odd.}
\end{cases}\\
z^-(\psi):=\begin{cases}
\mathrm{sc}(\psi),&\text{ if }\mathrm{sc}(\psi)\text{ is odd,}\\
\mathrm{sc}(\psi)+1,&\text{ if }\mathrm{sc}(\psi)\text{ is even.}
\end{cases}\label{eq:2.3}
\end{align}

The sign changes of a continuous function can indeed be infinitely many, however, under our assumptions this is not the case for periodic solutions of the DDE \eqref{eq:1.2}; see \cite[Theorem 2.4]{MPSe961}. From now, we drop the sign superscript in the notation of the zero number, our statements should then be read by picking $z^\pm$ when $f$ belongs to the corresponding subset $\mathfrak{X}^\pm$ of $\mathfrak{X}$. Let $x^\ast(t)$ and $x^\dagger(t)$ be two periodic solutions of the DDE \eqref{eq:1.2} with $f\in\mathfrak{X}$. Denote the minimal periods by $p$ and $p^\dagger$, and the periodic orbits by $\gamma^\ast=\gamma(x^\ast_0)$ and $\gamma^\dagger=\gamma(x^\dagger_0)$, respectively. By \cite[Theorem 2.2]{MPSe961}, we have that $z(x^\ast_t-x^\dagger_t)$ and $z(\dot{x}^\ast_t)$ are nonincreasing functions of $t\in\mathbb{R}$ and drop strictly at multiple zeros. Thanks to this strict dropping property, the \textit{planar projection}
\begin{align}
\begin{array}{rcl}
P:C&\to &\mathbb{R}^2\\
\psi&\mapsto& (\psi(0),\psi(-1)),
\end{array}
\end{align}
possesses special features when evaluated along periodic solutions of \eqref{eq:1.2}. We summarize the most relevant properties of $P$ in the following lemma.
\begin{lemma}\label{Lemma2.1}
Consider two periodic solutions $x^\ast(t)$ and $x^\dagger(t)$ of the DDE \eqref{eq:1.2} with orbits $\gamma^\ast = \gamma(x_0^\ast)$ and $\gamma^\dagger = \gamma(x_0^\dagger)$. Then the following statements hold:
\begin{enumerate}
\item[(i)]The planar projection $P\gamma^\ast$ is an embedding of $\gamma^\ast \subset C$ into $\mathbb{R}^2$. In particular, $P\gamma^\ast$ is a $C^1$ simple closed curve in $\mathbb{R}^2$.
\item[(ii)] The planar projections of two different orbits $\gamma^\ast$ and $\gamma^\dagger$ do not intersect, that is, they are nested and
\begin{align}\label{eq:2.7}
P\gamma^\ast\cap P\gamma^\dagger\neq \emptyset\quad \text{if and only if} \quad\gamma^\ast=\gamma^\dagger.
\end{align}
\item[(iii)] $x^\ast(t)$ is sinusoidal, that is, it moves monotonically in between its positive maximum and negative minimum values, and it reaches them exactly once over every minimal period.
\item[(iv)]$x^\ast(t)$ has the odd symmetry
\begin{align}\label{eq:2.8}
x^\ast\left(t-\frac{p}{2}\right) = -x^\ast(t).
\end{align}
\item[(v)]$P\gamma^\ast$ contains $(0,0)$ in its interior region.
\end{enumerate}
\end{lemma}
\begin{proof}
Part (i) follows from \cite[Theorem 2.1]{MPSe962}. 
Part (ii) is \cite[Lemma 5.7]{MPSe962} and Part (iii) is \cite[Theorem 7.1]{MPSe962}.
Finally, Part (iv) follows from \cite[Theorem 7.2]{MPSe962} and Part (v) is \cite[Corollary 7.4]{MPSe962} for the special case in which $x(t) \equiv 0$ is the only equilibrium of the DDE \eqref{eq:1.2}.
\end{proof}

\begin{figure}
\includegraphics[width=\textwidth]{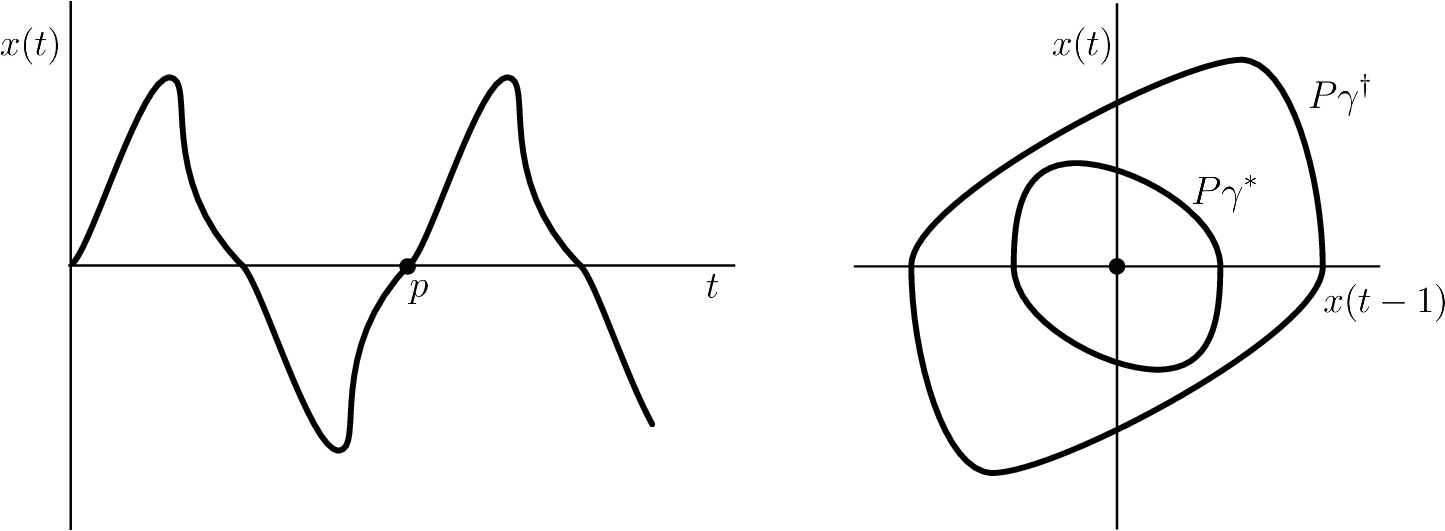}
\caption{(Left) Schematic picture for a sinusoidal periodic solution of the DDE \eqref{eq:1.2}. (Right) Nested planar projections of two different periodic solutions. }
\label{fig:2}
\end{figure}

Parts (i), (ii), (iv), and (v) in Lemma \ref{Lemma2.1} show that the planar projections of the periodic orbits are simple closed curves nested within one another. Moreover, they contain the origin $(0,0)$ in their interior regions and have a point symmetry with respect to the origin due to (iv). On the other hand, Part (iii) describes the shape of the solutions; see Figure \ref{fig:2}. 

We now turn our attention to the local stability analysis of the periodic solutions. Our presentation from this point follows \cite[Section 5]{MPNu13}. Recall that the local stability of the periodic orbit $\gamma^\ast$ comes characterized by the eigenvalues of the monodromy operator $M$, which coincides with the time-$p$ solution operator of the linearized equation \eqref{eq:1.10}.
Notice that the coefficients $\partial_1 f(x^\ast(t),x^\ast(t-1))$ and $\partial_2 f(x^\ast(t),x^\ast(t-1))$ have period $p/2$, rather than just $p$, due to the odd symmetry \eqref{eq:2.8} of $x^\ast(t)$ and the even-odd symmetry of $f\in\mathfrak{X}$. This motivates the study of the \textit{half-multipliers}, that is, the eigenvalues of the time-$(p/2)$ solution operator of \eqref{eq:1.10}, that we denote ${L}^{1/2}$. Observe that  $\dot{x}^\ast(t)$ is a periodic solution of \eqref{eq:1.10} and that by the odd symmetry \eqref{eq:2.8} we have $-1\in\mathrm{Spec}({L}^{1/2})$.

The zero number \eqref{eq:2.2}--\eqref{eq:2.3} induces a spectral ordering that relates the norm of the Floquet multipliers $\nu\in(0,\infty)$, to the zero number $z(\psi^\nu)$ of the associated real generalized eigenfunctions $\psi^\nu\in C$ of the operator ${L}^{1/2}$. This allows us to relate the unstable dimension $i(\gamma^\ast)$ to the zero number $z(\dot{x}^\ast_0)$ via the following lemma.

\begin{lemma}\label{Lemma2.2}
Let ${L}^{1/2}$ be the time-$(p/2)$ operator of a periodic solution $x^\ast(t)$ of the DDE \eqref{eq:1.2} with $f\in \mathfrak{X}$. There exists a real half-multiplier $-\mu_{\mathrm{c}}<0$ associated with 
\begin{align}
    \Psi_0\in \mathcal{E}_{\mu_\mathrm{c}} := \bigcup_{m\geq 1}\ker\left(-\mu_\mathrm{c} \mathrm{Id} - {L}^{1/2}\right)^m,
\end{align}
a real generalized eigenfunction satisfying $z(\Psi_0)=z(\dot{x}^\ast_0)$ and $\Psi_0\not\in \mathrm{span}_\mathbb{R}\{\dot{x}^\ast_0\}$. Moreover, exactly one of the following statements holds:
\begin{enumerate}
\item[(i)] $\mu_{\mathrm{c}}>1$, $\gamma^\ast$ is hyperbolic, $\dim \mathcal{E}_{\mu_\mathrm{c}}=1$, and $i(\gamma^\ast)=z(\dot{x}^\ast_0)$.
\item[(ii)] $\mu_{\mathrm{c}}< 1$, $\gamma^\ast$ is hyperbolic, $\dim \mathcal{E}_{\mu_\mathrm{c}}=1$, and $i(\gamma^\ast)=z(\dot{x}^\ast_0)-1$.
\item[(iii)] $\mu_{\mathrm{c}}= 1$, $\gamma^\ast$ is nonhyperbolic, $\dim \mathcal{E}_{\mu_\mathrm{c}}=2$, and $i(\gamma^\ast)=z(\dot{x}^\ast_0)-1$.
\end{enumerate}
\end{lemma}
\begin{proof}
The lemma follows from \cite[Theorem 5.1]{MPNu13}. We discuss the derivation briefly. Given $\nu\in (0,\infty)$, we define the direct sum of spaces
\begin{align}\label{eq:2.9}
\mathcal{G}_\nu:=\bigoplus_{|\mu|=\nu}\mathrm{Re}\left(\bigcup_{m \geq 1}\ker\left(\mu \mathrm{Id} - {L}^{1/2}\right)^m\right).
\end{align} 
Notice that, for real multipliers, $\mathcal{G}_{\nu}$ coincides with the associated eigenspace. Then \cite[Theorem 5.1]{MPNu13} shows that there exists a real multiplier $-\mu_{\mathrm{c}}<0$ such that $\dim\left(\mathcal{G}_{\mu_{\mathrm{c}}}+\mathcal{G}_{1}\right)=2$ and $z(\Psi_0)=z(\dot{x}^\ast_0)$ for all $\Psi_0\in \mathcal{G}_{\mu_{\mathrm{c}}} = \mathcal{E}_{\mu_{\mathrm{c}}}$. Moreover, counting multiplicities, the half-period monodromy operator ${L}^{1/2}$ possesses exactly $z(\dot{x}^\ast_0)-1$ eigenvalues $\mu$ such that  $|\mu|>\max\{1,\mu_{\mathrm{c}}\}$. Therefore, the unstable dimension $i(\gamma^\ast)$ satisfies
\begin{align}\label{eq:2.13}
i(\gamma^\ast)\in \{z(\dot{x}^\ast_0),z(\dot{x}^\ast_0)-1\}.
\end{align}

The case when $\dim \mathcal{G}_1=1$ corresponds to the periodic orbit $\gamma^\ast$ being hyperbolic with $\mathcal{G}_{1}=\mathrm{span}_{\mathbb{R}}\{\dot{x}_0^\ast\}$. Additionally, $\mu_{\mathrm{c}}\neq 1$ and $\dim\mathcal{G}_{\mu_{\mathrm{c}}}=1$. In virtue of \eqref{eq:2.13}, $i(\gamma^\ast)=z(\dot{x}^\ast_0)$ (resp., $i(\gamma^\ast)=z(\dot{x}^\ast_0)-1$) if $\mu_{\mathrm{c}}>1$ (resp., $\mu_{\mathrm{c}}<1$). This shows parts (i)--(ii).

The case $\dim \mathcal{G}_{1}=2$ characterizes a nonhyperbolic periodic orbit $\gamma^\ast$. Then we have $\mu_{\mathrm{c}}=1$ and $\mathcal{G}_{1}=\mathrm{span}_\mathbb{R}\{\dot{x}_0^\ast,\Psi_0\}$, where 
\begin{align}\label{eq:2.12}
\text{either} \quad \Psi_0\in\ker(\mathrm{Id} +{L}^{1/2}) \quad \text{or} \quad \Psi_0\in\ker(\mathrm{Id} + {L}^{1/2})^2\setminus\ker(\mathrm{Id} + {L}^{1/2}),
\end{align}
depending on the geometric multiplicity of $-1\in\mathrm{Spec}({L}^{1/2})$. Applying \eqref{eq:2.13} again, we obtain Part (iii).
\end{proof}

Therefore, the local stability of the periodic orbit $\gamma^\ast$ is completely determined by the half-multiplier $-\mu_{\mathrm{c}}$ in Lemma \ref{Lemma2.2}.  Moreover, any secondary bifurcation of periodic solutions from $x^\ast(t)$ corresponds to a change in the sign of $1-\mu_{\mathrm{c}}$. Thus, we say that $-\mu_{\mathrm{c}}$ is the \textit{critical half-multiplier} of $\gamma^\ast$.

\begin{remark}
    Lemma \ref{Lemma2.2} links the unstable dimension of $\gamma^\ast$ to $z(x^\ast_0)$. Since the solutions $x^\ast(t)$ are sinusoidal by Lemma \ref{Lemma2.1} (iii),  solutions with smaller minimal periods have a larger unstable dimension. Note that if $f\in\mathfrak{X}^+$, then $z(\dot{x}^\ast_0)\geq 2$ and the unstable dimension always satisfies $i(\gamma^\ast)\geq 1$. In particular, all periodic orbits of the DDE \eqref{eq:1.2} with $f\in\mathfrak{X}^+$ are unstable. However, the case $f\in\mathfrak{X}^-$ does allow the existence of asymptotically stable periodic orbits.
\end{remark}

\section{Proof of Theorem \ref{Theorem1.2}}\label{S3}

We begin the section by giving a proof of Lemma \ref{Lemma1.1}. Afterwards we prove the symmetries of the solutions of the ODEs \eqref{eq:1.12} in Lemma \ref{Lemma3.1}. Lemma \ref{Lemma3.2} relates the  period of a periodic orbit $\gamma^\ast$ of the DDE \eqref{eq:1.2} with $f\in\mathfrak{X}$ to the geometrical shape of its planar projection $P\gamma^\ast$. We conclude the section with the proof of Theorem \ref{Theorem1.2}.

\begin{proof}[Proof of Lemma \ref{Lemma1.1}]
Since $f\in\mathfrak{X}$, notice that $f({u_1},{u_2})=0$ if and only if ${u_2}=0$, for all ${u_1}\in\mathbb{R}$. Therefore, $f({u_1},0)=f(0,{u_1})=0$ if and only if ${u_1}=0$ and $(0,0)$ is the only equilibrium of \eqref{eq:1.12}.

Next we show that the solution $({u_1}(t; a),{u_2}(t; a))$ of the ODEs \eqref{eq:1.12} with the initial condition $({u_1}(0; a),{u_2}(0; a))=(a, 0)$, $a>0$, is periodic. We recall that the vector field $F({u_1},{u_2}):=(f({u_1},{u_2}),-f({u_1},{u_2}))$ possesses the reversibility $\sigma$ defined in \eqref{eq:1.8}. Linearizing \eqref{eq:1.12} at $(0,0)$ we observe that $\mathrm{Spec}(D_{(0,0)}F)=\{i,-i\}$. We apply a reversible Hopf bifurcation theorem at $(0,0)$; see \cite[Theorem 7.5.4.]{Van82}. Hence, there exists a continuum of periodic solution surrounding the origin. Since $(0,0)$ is the only equilibrium of \eqref{eq:1.12}, if the continuum is bounded, then, by the planar Poincar\'{e}--Bendixson theorem, it has a periodic solution at the boundary with gradient dynamics in the exterior region. However, this scenario is in contradiction with the reversibilities \eqref{eq:1.8} and shows that the family of periodic solutions covers $\mathbb{R}^2$.

In particular, the map ${p}_f$ assigning the minimal period ${p}_f(a)$ to the periodic solution $({u_1}(t; a),{u_2}(t; a))$ is well defined. Furthermore, since $\dot{u}_2({p}_f(a); a)=-f(0,a)\neq 0$ for $a> 0$, ${p}_f$ coincides with the solution obtained by using the implicit function theorem to solve
\begin{align}
{u_2}({p}_f(a); a)-{u_2}(0; a)=0.
\end{align}
In particular, ${p}_f$ inherits its regularity from $f$.
\end{proof}

\begin{lemma}\label{Lemma3.1}
Consider the planar ODEs \eqref{eq:1.12} with $f\in\mathfrak{X}$. Let $({u_1}(t; a),{u_2}(t; a))$ denote the solution with the initial condition $(a, 0)$, $a>0$. Then 
\begin{align}
({u_1}(t; a),{u_2}(t; a))&=\left(-{u_2}\left(t-\frac{3{p}_f(a)}{4}; a\right), {u_1}\left(t-\frac{3{p}_f(a)}{4}; a\right)\right)\;\text{if}\quad f\in\mathfrak{X}^+,\label{eq:3.2}\\
({u_1}(t; a),{u_2}(t; a))&=\left(-{u_2}\left(t-\frac{{p}_f(a)}{4}; a\right), {u_1}\left(t-\frac{{p}_f(a)}{4}; a\right)\right)\;\text{if}\quad f\in\mathfrak{X}^-.\label{eq:3.3}
\end{align}
Moreover, ${u_1}(t; a)$ satisfies the symmetry
\begin{align}\label{eq:3.3b}
{u_1}\left(t-\frac{{p}_f(a)}{2}; a\right)=-{u_1}(t; a).
\end{align}
\end{lemma}

\begin{proof}
By Lemma \ref{Lemma1.1}, the orbits $O_a:= \{({u_1}(t; a),{u_2}(t; a)) : t\in\mathbb{R}\}\subset \mathbb{R}^2$ are simple curves surrounding $(0,0)$. By the equivariance of the ODEs \eqref{eq:1.12} under the rotation $\varrho$ defined in \eqref{eq:1.7}, we have that $O_a$ are $\varrho$-invariant, that is, $\varrho (O_a)=O_a$. 

Since $\varrho({u_1}(t; a),{u_2}(t; a))=({u_2}(t; a),-{u_1}(t; a))$ solves the ODEs \eqref{eq:1.12} and shares orbits with $({u_1}(t; a),{u_2}(t; a))$, we have that ${u_2}(t; a)={u_1}(t-\tau; a)$ for some $\tau\in [0,{p}_f(a)]$. Moreover, we have that $4\tau$ is a multiple of the period because $\varrho^4=\mathrm{Id}$. Therefore, $\tau$ can only take the values ${p}_f(a)/4$ and $3{p}_f(a)/4$. In the case $f\in\mathfrak{X}^+$, the points in $O_a$ wind in a clockwise direction around the origin $(0,0)$ as $t$ grows. Thus ${u_2}(t; a)$ decreases until it reaches its minimum at $-a$ and then grows until it reaches its maximum $a$. We conclude that ${u_1}(0; a)={u_2}(3{p}_f(a)/4; a)=a$ and $\tau=3{p}_f(a)/4$. If $f\in\mathfrak{X}^-$, then the winding is clockwise and a similar argument yields $\tau={p}_f(a)/4$. This proves the expressions \eqref{eq:3.2}--\eqref{eq:3.3}. Using the periodicity of $({u_1}(t; a),{u_2}(t; a))$ and applying \eqref{eq:3.2}--\eqref{eq:3.3} twice, we obtain \eqref{eq:3.3b}.
\end{proof}

\begin{lemma}\label{Lemma3.2}
Let $x^\ast(t)$ be a periodic solution of the DDE \eqref{eq:1.2} with nonlinearity $f\in\mathfrak{X}$ and let $\gamma^\ast = \gamma^\ast(x^\ast_0) \subset C$ denote the corresponding orbit. Then $x^\ast(t)$ has period $4$ if and only if the planar projection $P\gamma^\ast=\{(x^\ast(t),x^\ast(t-1)) : t\in\mathbb{R}\}$ intersects the vertical axis orthogonally.
\end{lemma}
\begin{proof}
First, by \cite[Lemma 4.1]{ChMP78}, we know that $2$ is not a period of $x^\ast(t)$. Indeed, otherwise the ODEs
\begin{align}\label{eq:3.6a}
\begin{split}
\dot{u}_1&=f({u_1},{u_2}),\\
\dot{u}_2&=f({u_2},{u_1}),
\end{split}
\end{align}
would possess a periodic solution $(x^\ast(t),x^\ast(t-1))$. However, this is impossible since the diagonal $\{({u_1},{u_1}):{u_1}\in\mathbb{R}\}$ is invariant under the dynamics of the ODEs \eqref{eq:3.6a}.

Suppose that $x^\ast(t)$ has minimal period $p$ and period $4$. Since $2$ is not a period of $x^\ast(t)$, we have that $(2m-1)p=4$ for some $m\in\mathbb{N}$. Then $x^\ast(t-2)=x^\ast(t-mp+mp/2)=-x^\ast(t)$ by the odd symmetry \eqref{eq:2.8} and $(x^\ast(t),x^\ast(t-1))$ is a solution of the planar ODEs \eqref{eq:1.12}. In particular, $P\gamma^\ast$ is an orbit of the ODEs \eqref{eq:1.2} and, since $f\in\mathfrak{X}$, evaluating \eqref{eq:1.12} at the intersections with the vertical axis shows that they must happen orthogonally.

To see the converse implication, assume that the simple curve $P\gamma^\ast$ intersects the vertical axis orthogonally. We set without loss of generality $x^\ast(0)=\bar{x}:= \max_{t\in \mathbb{R}}x^\ast(t)$ and claim that $x^\ast(-2)=-x^\ast(0)=-\bar{x}$.
Indeed, $x^\ast(0)=\bar{x}$ is a maximum, and thus
\begin{align}\label{eq:3.6}
\dot{x}^\ast(0)=f(x^\ast(0),x^\ast(-1))=0.
\end{align}
For $f\in\mathfrak{X}$, \eqref{eq:3.6} implies that $x^\ast(-1)=0$. Lemma \ref{Lemma2.1} (iii) states that $x^\ast(t)$ moves monotonically between extrema. Therefore, the point $(x^\ast(-1),x^\ast(-2))=(0,x^\ast(-2))$ corresponds to one of the two intersections that $P\gamma^\ast$ has with the vertical axis. As a result,
\begin{align}
\dot{x}^\ast(-2)=f(x^\ast(-2),x^\ast(-3))=0,
\end{align}
and $f\in\mathfrak{X}$ implies $x^\ast(-3)=0$. By Lemma \ref{Lemma2.1} (iii) and the odd symmetry \eqref{eq:2.8}, $P\gamma^\ast$ intersects the horizontal axis at $(\bar{x},0)$ and $(-\bar{x},0)$, only. Since $2$ is not a period for $x^\ast(t)$, the only possibility is $x^\ast(-2)=-\bar{x}=-x^\ast(0)$. Repeating the procedure once more shows that $x^\ast(-4)=x^\ast(0)$ and completes the proof.
\end{proof}

\begin{proof}[Proof of Theorem \ref{Theorem1.2}]
Let $x^\ast(t)$ be a solution of the DDE \eqref{eq:1.2} with $f\in\mathfrak{X}$ and minimal period $p$. We denote by $\gamma^\ast = \gamma^\ast(x_0^\ast)\subset C$ the orbit of $x^\ast(t)$. Notice that $x^\ast(t)$ solves the family of DDEs $\dot{x}^\ast(t)=f(x^\ast(t),x^\ast(t-1-mp))$ for all $m\in\mathbb{Z}$. Rescaling time, we have that $x^{(m)}(t):=x^\ast((1+mp)t)$ solves the DDE
\begin{align}\label{eq:3.8}
\dot{x}(t)=(1+mp)f(x(t),x(t+1)).
\end{align}
Moreover, for all $m\in \mathbb{Z}$, the planar projections $P\gamma(x^{(m)}_0)$ satisfy
\begin{align}\label{eq:3.9}
\begin{split}
P\gamma\left(x^{(m)}_0\right)&= \{(x^{(m)}(t),x^{(m)}(t-1)) : t\in\mathbb{R}\}\\
&= \{(x^\ast(t),x^\ast(t-1-mp)) : t\in\mathbb{R}\}\\
&=\{(x^\ast(t),x^\ast(t-1) : t\in\mathbb{R}\}\\
&=P\gamma^\ast.
\end{split}
\end{align}
Since the term $1+mp$ in the DDE \eqref{eq:3.8} can be chosen negative, we restrict our attention to $f\in\mathfrak{X}^+$ without loss of generality. The feedback type will only play a role towards the end of the proof, when we prove the identities \eqref{eq:1.14}--\eqref{eq:1.15}.

We proceed by contradiction.  Suppose that $x^\ast(t)$ does not have period $4$, that is, that $4$ is not a multiple of $p$. By Lemma \ref{Lemma3.2}, the intersection of the projection $P\gamma^\ast=\{(x^\ast(t),x^\ast(t-1)) : t\in\mathbb{R}\}$ with the vertical axis is not orthogonal. We denote the intersection point
\begin{align}\label{eq:3.10}
(0,\hat{x})=P\gamma^\ast\cap \{(0,{u_2}) : {u_2}\in\mathbb{R}^+\}.
\end{align}
Let us denote by $({u_1}(t; a),{u_2}(t; a))$ the periodic solution of the ODEs \eqref{eq:1.12} with the initial condition $(a, 0)$, $a\geq 0$. Since $f\in\mathfrak{X}$, the orbit
\begin{align}\label{eq:3.11}
O_a:=\{({u_1}(t; a),{u_2}(t; a)) : t\in\mathbb{R}\}
\end{align}
intersects both the horizontal and the vertical axes orthogonally. Therefore, we can find an $\varepsilon>0$ such that 
\begin{align}\label{eq:3.12}
P\gamma^\ast \cap O_a \neq \emptyset, \quad\text{for all } a\in(\hat{x}-\varepsilon,\hat{x}+\varepsilon).
\end{align}
By Lemma \ref{Lemma3.1} with $f\in\mathfrak{X}^+$, we know that ${u_1}(t; a)$ satisfies the DDE
\begin{align}\label{eq:3.13}
\dot{u}_1(t; a)&=f\left({u_1}(t; a),{u_1}\left(t-\frac{3{p}_f(a)}{4}; a\right)\right).
\end{align}
If $3{p}_f(a_0)/4=1$ for some $a_0\in(\hat{x}-\varepsilon,\hat{x}+\varepsilon)$, the proof is finished because $O_{a_0}$ corresponds to the planar projection of the periodic solution ${u_1}(t; a_0)$ of the DDE \eqref{eq:1.2} and $O_{a_0}$ intersects $P\gamma^\ast$, which is a contradiction to the nesting property in Lemma \ref{Lemma2.1} (ii).

We now discuss the case $3{p}_f(a)/4\neq 1$ for all $a\in (\hat{x}-\varepsilon,\hat{x}+\varepsilon)$. By the argument at the beginning of the proof, for all $n \in \mathbb{N}$, the functions 
\begin{align}
{u}^{(a,n)}_1(t):={u_1}\left(\left(n - \frac{1}{4}\right){p}_f(a)t; a\right),
\end{align}
solve the DDE
\begin{align}\label{eq:3.15}
\dot{x}(t)={r} f(x(t), x(t-1))
\end{align}
for parameter values 
\begin{align}
{r}=\left(n - \frac{1}{4}\right){p}_f(a).
\end{align} 
Proceeding as in \eqref{eq:3.9}, we denote $\psi^{(n)}(\theta) := {u}^{(a,n)}_1(\theta)$ for $\theta\in[-1, 0]$ and see that
\begin{align}\label{eq:3.16}
P\left(\psi^{(n)}\right) = O_{a}, \quad \text{for all}\; n\in\mathbb{N}.
\end{align}

Our argument at this point becomes very geometrical, the main idea is that the branches  $\mathcal{B}^+_n:=\{({r},a) :  {r}=(n - 1/4){p}_f(a)\}$ of periodic solutions of the DDE \eqref{eq:3.15} become steeper and
cover larger ${r}$-regions for larger $n\in\mathbb{N}$; see Figure \ref{fig:3}.

 Recall that the period map ${p}_f(a)$ is assumed to be locally nonconstant. Thus, for all $\delta > 0$ we can find $m^\ast, n^\ast \in\mathbb{N}$ such that
\begin{align}\label{eq:3.17}
\left\lvert\left(n^\ast - \frac{1}{4}\right){p}_f(\hat{x}) - 1 - m^\ast {p} \right\rvert < \delta.
\end{align}
If $\delta>0$ is small enough, then ${p}_f$ being locally nonconstant implies that there exists an $a_0\in(\hat{x}-\varepsilon,\hat{x}+\varepsilon)$
\begin{align}
  \left(n^\ast - \frac{1}{4}\right){p}_f(a_0) = 1 + m^\ast p.
\end{align}
Therefore, we have that $x^{(m^\ast)}(t)$ and ${u}^{(a_0,n^\ast)}_1(t)$ are both solutions of \eqref{eq:3.15} with the parameter value ${r}=1 + m^\ast p = (n^\ast - 1/4)p_f(a_0)$. By the identities \eqref{eq:3.9}, \eqref{eq:3.12}, and \eqref{eq:3.16}, we have that $P\gamma(x^{(m^\ast)}_0)\cap P\gamma(\psi^{(n^\ast)})=P\gamma^\ast\cap O_{a_0}\neq \emptyset$, in contradiction to Lemma \ref{Lemma3.1} (ii). Therefore, $x^\ast(t)$ has period $4$.

\begin{figure}
\includegraphics[width=\textwidth]{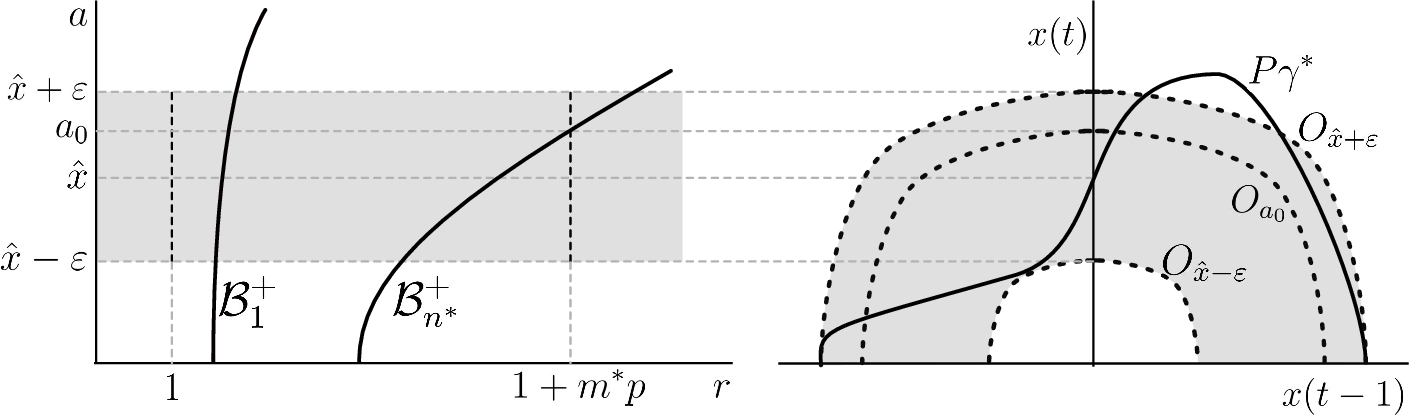}
\caption{(Left) Branches of periodic solutions $\mathcal{B}^+_n$ of the DDE \eqref{eq:3.15} (solid). At parameter value $r = 1 + m^\ast p$, there are periodic solutions on the branch $\mathcal{B}^+_{n^\ast}$ whose projection intersects $P\gamma^\ast$ in the $(x(t), x(t - 1))$-plane. 
(Right) Planar projection $P\gamma^\ast$ of the orbits $\gamma(x_0^{(m)})$ for all $m\in\mathbb{Z}$ (solid) and orbits of the ODEs \eqref{eq:1.12} (dashed). All the periodic orbits of \eqref{eq:1.12} possessing amplitudes within the gray region intersect the planar projection $P\gamma^\ast$. }
\label{fig:3}
\end{figure}

Having period $4$ implies that $(x^\ast(t),x^\ast(t-1))$ solves the ODE \eqref{eq:1.12} with minimal period ${p}_f(\bar{x})$, where $\bar{x}=\max_{t\in[0,4]}{x^\ast(t)}$. Moreover, by Lemma \ref{Lemma3.1}, we know that the minimal periods satisfy
\begin{align}
1=\left(n - \frac{1}{4}\right){p}_f(\bar{x}),\quad \text{for some }n\in\mathbb{N},
\end{align}
which yields \eqref{eq:1.14}. For the case $f\in\mathfrak{X}^-$, Lemma \ref{Lemma3.1} yields
\begin{align}
1=\left(n - \frac{3}{4}\right){p}_f(\bar{x}),\quad \text{for some }n\in\mathbb{N},
\end{align}
which proves \eqref{eq:1.15}. In showing \eqref{eq:1.14}--\eqref{eq:1.15}, we have also proved the converse part of Theorem \ref{Theorem1.2}, concluding the proof.
\end{proof}

\section{Proof of Theorem \ref{Theorem1.3}}\label{S4}
By Lemma \ref{Lemma2.2}, we have to discuss the sign of $1-\mu_{\mathrm{c}}$ where $-\mu_{\mathrm{c}}$ is the critical half-multiplier of the periodic orbit $\gamma^\ast$, as defined in Lemma \ref{Lemma2.2}. We begin by establishing a connection between the unstable dimension $i(\gamma^\ast)$ and the minimal period of the periodic solution in Lemma \ref{Lemma4.1} and Corollary \ref{Corollary4.2}. Then, we characterize the hyperbolicity of the periodic orbit in Lemma \ref{Lemma4.3}. Finally, we prove Theorem \ref{Theorem1.3}.

\begin{lemma}\label{Lemma4.1}
Consider a periodic solution $x^\ast(t)$ of the DDE \eqref{eq:1.2} with nonlinearity $f\in\mathfrak{X}$. If $f\in \mathfrak{X}^+$ and $x^\ast(t)$ has minimal period $p:= 4/(4n-1)$ for some $n\in \mathbb{N}$, then $z(\dot{x}^\ast_0)=2n$. If $f\in \mathfrak{X}^-$ and $x^\ast(t)$ has minimal period $p:= 4/(4n-3)$ for some $n\in \mathbb{N}$, then $z(\dot{x}^\ast_0)=2n-1$.
\end{lemma}
\begin{proof}
We know from Lemma \ref{Lemma2.1} (iii) that $x^\ast(t)$ attains its maximum $\bar{x}$ (resp., its minimum $\underline{x}$) once over every minimal period and moves monotonically between the maximum and the minimum. By the odd symmetry \eqref{eq:2.8}, any two neighboring sign changes of $\dot{x}^\ast(t)$ are separated by half the minimal period. Therefore, if $f\in\mathfrak{X}^+$, then due to $p=4/(4n-1)$ we have the bounds
\begin{align*}
2n-1= \left\lfloor\frac{2}{p}\right\rfloor \leq  z(\dot{x}^\ast_0)\leq\left\lceil\frac{2}{p}\right\rceil=2n,
\end{align*}
and if $f\in\mathfrak{X}^-$, then due to $p=4/(4n-3)$,
\begin{align*}
2n-2 = \left\lfloor\frac{2}{p}\right\rfloor \leq  z(\dot{x}^\ast_0) \leq\left\lceil\frac{2}{p}\right\rceil = 2n-1.
\end{align*}
Here we used the notation $\left\lfloor\cdot\right\rfloor$ (resp., $\left\lceil\cdot\right\rceil$) for the standard floor (resp., ceiling) function. The parity choice in the definition of the zero number \eqref{eq:2.2}--\eqref{eq:2.3} shows that $z(x^\ast_0)=2n$ (resp., $z(x^\ast_0)=2n-1$) if $f\in\mathfrak{X^+}$ (resp., $f\in\mathfrak{X}^-$).
\end{proof}

\begin{corollary}\label{Corollary4.2}
In the setting of Lemma \ref{Lemma4.1}, let $\gamma^\ast = \gamma(x_0^\ast) \subset C$ be the orbit of $x^\ast(t)$. We denote by $-\mu_{\mathrm{c}}$ the critical half-multiplier of $\gamma^\ast$. Then exactly one of the following statements holds:
\begin{itemize}
\item[(i)]$\mu_{\mathrm{c}}>1$, $\gamma^\ast$ is hyperbolic, and $i(\gamma^\ast)= 2n$ (resp., $i(\gamma^\ast)=2n-1$) if $f\in\mathfrak{X^+}$  (resp., $f\in\mathfrak{X^-}$).
\item[(ii)]$\mu_{\mathrm{c}}<1$, $\gamma^\ast$ is hyperbolic, and $i(\gamma^\ast)= 2n-1$ (resp., $i(\gamma^\ast)=2n-2$) if $f\in\mathfrak{X^+}$  (resp., $f\in\mathfrak{X^-}$).
\item[(iii)]$\mu_{\mathrm{c}}=1$, $\gamma^\ast$ is nonhyperbolic, and $i(\gamma^\ast)= 2n-1$ (resp., $i(\gamma^\ast)=2n-2$) if $f\in\mathfrak{X^+}$  (resp., $f\in\mathfrak{X^-}$).
\end{itemize}
\end{corollary}
\begin{proof}
We plug Lemma \ref{Lemma4.1} into the three cases we distinguished in Lemma \ref{Lemma2.2}.
\end{proof}

\begin{lemma}\label{Lemma4.3}
In the setting of Lemma \ref{Lemma4.1}, let $\bar{x}:= \max_{t\in\mathbb{R}}x^\ast(t)$ and denote by $\gamma^\ast$ the orbit of $x^\ast(t)$. Then $\gamma^\ast$ is hyperbolic if and only if ${p}_f'(\bar{x})\neq 0$.
\end{lemma}

\begin{proof}
The plan for the proof is the following.
To see how ${p}_f'(\bar{x})=0$ implies that $\gamma^\ast$ is nonhyperbolic we construct a periodic solution of the linearized equation \eqref{eq:1.10} explicitly. We do this in such a way that the newly obtained solution is linearly independent from $\dot{x}^\ast(t)$. To show the converse, we first prove that if $\gamma^\ast$ is not hyperbolic, then any generalized eigenfunction is genuinely an eigenfunction, that is, the trivial eigenvalue $1$ of the monodromy operator ${L}$ defined in \eqref{eq:1.9} has geometric multiplicity two. In particular, this fact implies that the period map must have a critical point, by examining an associated boundary value problem.

Denoting ${p}:= {p}_f(\bar{x})$ and ${p}':= {p}'_f(\bar{x})$, we first suppose that ${p}'=0$. Let us denote by $x^\ast(t; a)$, the solution of the DDE
\begin{align}\label{eq:4.1}
\dot{x}(t)=f\left(x(t),x\left(t-\frac{m{p}_f(a)}{4}\right)\right),
\end{align} 
with minimal period ${p}_f(a)$ and amplitude $a= \max_{t\in\mathbb{R}}x^\ast(t; a)$. Here $m \in \mathbb{N}$ is such that $m p / 4 = 1$ and we have normalized the initial condition so that 
\begin{align}\label{eq:4.2}
x^\ast(0; a)=a.
\end{align}
We may assume, shifting time if necessary, that $x^\ast(t)=x^\ast(t; \bar{x})$. Then the amplitude derivative $y^\ast(t; a):= \partial_a x^\ast(t; a)$ solves the linear inhomogeneous nonautonomous DDE
\begin{align}\label{eq:4.3}
\begin{split}
\dot{y}(t) &= \partial_1f\left(x^\ast(t; a),x^\ast\left(t-\frac{m{p}_f(a)}{4}; a\right)\right)y(t)\\
&\quad+\partial_2f\left(x^\ast(t; a),x^\ast\left(t-\frac{m{p}_f(a)}{4}; a\right)\right)y\left(t-\frac{m{p}_f(a)}{4}\right)\\
&\quad- \partial_2f\left(x^\ast(t; a),x^\ast\left(t-\frac{m{p}_f(a)}{4}; a\right)\right) \partial_t {x}^\ast\left(t-\frac{m{p}_f(a)}{4}; a\right)\frac{m{p}_f'(a)}{4}.
\end{split}
\end{align}
Since ${p}'=0$, then $y^\ast(t):= y^\ast(t; \bar{x})$ satisfies the linearized equation \eqref{eq:1.10} around $x^\ast(t)$. In addition, $y^\ast(t)$ is periodic since
\begin{align}\label{eq:4.4}
\begin{split}
y^\ast(t+{p})&=\partial_a x^\ast(t+{p}; \bar{x})+{p}'\partial_t x^\ast(t + {p}; \bar{x}))\\
&=\partial_a x^\ast(t+{p}; \bar{x})\\
&=y^\ast(t).
\end{split}
\end{align}
Clearly, $y^\ast(t)$ satisfies $y^\ast(0)=1$, in contrast to the trivial solution of the linearized equation $\dot{x}^\ast(t)$ for which $\dot{x}^\ast(0)=0$ due to the normalization \eqref{eq:4.2}. Therefore, $\dot{x}^\ast_0$ and $y^\ast_0$ are linearly independent. In other words, the Floquet multiplier $1$ of $\gamma^\ast$ has geometric multiplicity at least $2$ and so $\gamma^\ast$ is not hyperbolic.

Suppose that $\gamma^\ast$ is not hyperbolic. By Lemma \ref{Lemma2.2} (iii), the half-multiplier $-1\in\mathrm{Spec}({L}^{1/2})$ has algebraic multiplicity two. Here ${L}^{1/2}$ is the time-$({p}/2)$ solution operator of the linearized equation \eqref{eq:1.10}. 

Let us denote by $\Psi(t)$, the solution of the linearized equation \eqref{eq:1.10} with initial condition the eigenfunction $\Psi_0$ associated to $-\mu_\mathrm{c}$. By the Floquet theory in \cite{HaLu93}, $\Psi(t)$ has the form
\begin{align}\label{eq:4.5}
\Psi(t)=v^\ast(t)+\kappa t\dot{x}^\ast(t),
\end{align}
where $v^\ast(t)$ has period ${p}$ and satisfies $v^\ast(t-{p}/2)=-v^\ast(t)$ for $\kappa\in\mathbb{R}$. Note that $m{p}/2=2$ and therefore $v^\ast(t-2)=-v^\ast(t)$. In the following we denote $A(t):= \partial_1 f(x^\ast(t),x^\ast(t-1))$ and $B(t):= \partial_2 f(x^\ast(t),x^\ast(t-1))$.
In particular, $(v^\ast(t),v^\ast(t-1))$ solves the inhomogeneous ODEs
\begin{align}\label{eq:4.6}
\begin{split}
\dot{v}_1(t)&=A(t)v_1(t)+B(t)v_2(t)-\kappa B(t)\dot{x}^\ast(t-1)-\kappa\dot{x}^\ast(t),\\
\dot{v}_2(t)&=-B(t-1)v_1(t)+A(t-1)v_2(t)+\kappa B(t-1)\dot{x}^\ast(t)-\kappa \dot{x}^\ast(t-1).
\end{split}
\end{align}
Let us use the transpose superscript to highlight the solutions of the adjoint equation of the ODEs \eqref{eq:4.6} with respect to the Euclidean inner product. By a Fredholm alternative argument in \cite[Section IV.1, Lemma 1.1]{Ha69}, \eqref{eq:4.6} have solutions with period ${p}$ if and only if
\begin{align}\label{eq:4.7}
\kappa\int_0^{p} (v_1^{\mathsf{T}}(t),v_2^{\mathsf{T}}(t))
\left(
\begin{array}{c}
-B(t)\dot{x}^\ast(t-1)-\dot{x}^\ast(t)\\
B(t-1)\dot{x}^\ast(t)-\dot{x}^\ast(t-1)
\end{array}\right)\mathrm{d}t=0.
\end{align}
For all $(v_1^{\mathsf{T}}(t),v_2^{\mathsf{T}}(t))$ with period ${p}$ solving the adjoint equation
\begin{align}\label{eq:4.8}
\begin{split}
\dot{v}_1^{\mathsf{T}}(t)&=-A(t)v_1^{\mathsf{T}}(t)+B(t-1)v_2^{\mathsf{T}}(t),\\
\dot{v}_2^{\mathsf{T}}(t)&=-B(t)v_1^{\mathsf{T}}(t)-A(t-1)v_2^{\mathsf{T}}(t).
\end{split}
\end{align}
In particular, we have that
\begin{align}\label{eq:4.9}
(v_1^{\mathsf{T}}(t),v_2^{\mathsf{T}}(t)):= \mathrm{e}^{-\int_{t-1}^tA(s)\mathrm{d}s}\left(-\dot{x}^\ast(t-1)
, \dot{x}^\ast(t)
\right),
\end{align}
is a ${p}$-periodic solution of the adjoint equation \eqref{eq:4.8}. Plugging the solution \eqref{eq:4.9} into the condition \eqref{eq:4.7}, we obtain that
\begin{align}\label{eq:4.10}
\kappa \int_0^{p} \mathrm{e}^{-\int_{t-1}^t A(s)\mathrm{d}s}\left(B(t)|\dot{x}^\ast(t-1)|^2+B(t-1)|\dot{x}^\ast(t)|^2\right)\mathrm{d}t=0.
\end{align}
 By the monotone feedback condition $f\in\mathfrak{X}$, we have that $B(t)\neq 0$ for all $t\in\mathbb{R}$; thus the integral term in \eqref{eq:4.10} never vanishes and we conclude that $\kappa=0$.

Hence, the generalized eigenfunction $\Psi(t)= v^\ast(t)$ is a periodic solution of the linearized equation \eqref{eq:1.10} and $\Psi_0$ is linearly independent of $\dot{x}_0^\ast$, by construction. Moreover, $\Psi(t-2)=-\Psi(t)$, which implies that the linear ODEs
\begin{align}\label{eq:4.11}
\begin{split}
\dot{v}_1(t)&=A(t)v_1(t)+B(t)v_2(t),\\
\dot{v}_2(t)&=-B(t-1)v_1(t)+A(t-1)v_2(t),
\end{split}
\end{align}
have two linearly independent solutions, namely $(\dot{x}^\ast(t),\dot{x}^\ast(t-1))$ and $(\Psi(t),\Psi(t-1))$, with period $4$.

For $a\geq 0$, we denote by $({u_1}(t; a),{u_2}(t; a))$ the solution of the ODEs \eqref{eq:1.12} with the initial condition $(a, 0)$, that is,
\begin{align}\label{eq:4.12}
\begin{split}
\dot{u}_1(t; a)&=f({u_1}(t; a),{u_2}(t; a)),\\
\dot{u}_2(t; a)&=-f({u_2}(t; a),{u_1}(t; a)),
\end{split}
\end{align}
satisfying the boundary condition
\begin{align}\label{eq:4.13}
(a, 0)=({u_1}(0; a),{u_2}(0; a))=({u_1}({p}_f(a); a), {u_2}({p}_f(a); a)).
\end{align}
By oddness of the periodic solution \eqref{eq:2.8} and even-odd symmetry of $f\in\mathfrak{X}$, we have that the amplitude derivative $(w_1(t),w_2(t)):= \partial_a({u_1}(\bar{x},t),{u_2}(\bar{x},t))$ solves the linear ODEs \eqref{eq:4.11}. Moreover, it takes the boundary values
\begin{align}\label{eq:4.14}
\begin{split}
(1,0)&=(w_1(0),w_2(0))\\
&=(w_1(4),w_2(4))+m{p}'(\dot{x}^\ast(0),\dot{x}^\ast(-1)).
\end{split}
\end{align}
However, $(w_1(t),w_2(t))$ can be written as a linear combination of the eigenfunctions $(\dot{x}^\ast(t),\dot{x}^\ast(t-1))$ and $(\Psi(t),\Psi(t-1))$ of \eqref{eq:4.11}, both of which are ${p}$-periodic and also $4$-periodic. Since $\dot{x}^\ast(-1)\neq 0$, by the normalized form \eqref{eq:4.2}, it follows that ${p}'=0$ so that \eqref{eq:4.14} is satisfied.
\end{proof}

\begin{proof}[Proof of Theorem \ref{Theorem1.3}]
We use the notation ${p}:= {p}_f(\bar{x})$ and ${p}':= {p}'_f(\bar{x})$. We first deal with the result in the positive feedback case, that is, let $x^\ast(t)$ be a periodic solution with orbit $\gamma^\ast$ of the DDE \eqref{eq:1.2} with nonlinearity $f\in\mathfrak{X^+}$ and period $4$. In particular, $(x^\ast(t),x^\ast(t-1))$ solves the ODEs \eqref{eq:1.12} and has minimal period ${p}=4/(4n-1)$ for some $n\in\mathbb{N}$. Without loss of generality, we normalize $x^\ast(t)$ so that $x^\ast(0)=\bar{x}$.

By Lemma \ref{Lemma4.3}, if ${p}'= 0$, then the critical half-multiplier $-\mu_{\mathrm{c}}$ of the periodic orbit $\gamma^\ast$ satisfies $-\mu_{\mathrm{c}}=-1$ and the time-$({p}/2)$ solution operator ${L}^{1/2}$ has a geometrically double eigenvalue $-1$. By Corollary \ref{Corollary4.2} (iii), the unstable dimension of $\gamma^\ast$ is 
\begin{eqnarray}\label{eq:4.15}
i(\gamma^\ast)=2n-1.
\end{eqnarray}
If ${p}'\neq0$, we know by Lemma \ref{Lemma4.3} that $\gamma^\ast$ is hyperbolic and Lemma \ref{Lemma2.2} (i)--(ii) imply that $-\mu_{\mathrm{c}}<0$ has geometric multiplicity one. Furthermore, the associated eigenfunction $\Psi_0$ satisfies $z(\Psi_0)\equiv z(\dot{x}^\ast_0)$. 

In virtue of Corollary \ref{Corollary4.2}, we have that $i(\gamma^\ast)=2n$ (resp., $i(\gamma^\ast)=2n-1$) if $\mu_{\mathrm{c}}>1$ (resp., $\mu_{\mathrm{c}}<1$).
We now show by a comparison argument that if ${p}'>0$, then $\mu_{\mathrm{c}}<1$. 

Let $\Psi(t)$ denote the solution of the linearized equation \eqref{eq:1.10} whose initial condition is the eigenfunction $\Psi_0$ associated to the critical half-multiplier $-\mu_{\mathrm{c}}$. Then, we have that
\begin{align}\label{eq:4.16}
\Psi\left(t-\frac{{p}}{2}\right)=-\frac{1}{\mu_\mathrm{c}}\Psi(t),
\end{align}
recalling that $2=(4n-1){p}/2$, we define $ \nu:= (\mu_{\mathrm{c}})^{(4n-1){p}/2}>1$. This choice yields $\Psi(t-2)=-\nu^{-1}\Psi(t)$ and 
\begin{align}
\begin{split}
\mu_{\mathrm{c}}>1& \quad \text{if and only if} \quad \nu>1,\\
\mu_{\mathrm{c}}<1& \quad \text{if and only if} \quad \nu<1.
\end{split}
\end{align}
In particular, $(v_1(t),v_2(t))=(\Psi(t),\Psi(t-1))$ solves the ODEs
\begin{align}\label{eq:4.18}
\begin{split}
\dot{v}_1(t)&=A(t)v_1(t)+B(t)v_2(t),\\
\dot{v}_2(t)&=-\frac{1}{\beta}B(t-1)v_1(t)+A(t-1)v_2(t),
\end{split}
\end{align}
for the parameter value $\beta=\nu$; here we used the notation $A(t):= \partial_1 f(x^\ast(t),x^\ast(t-1))$ and $B(t):= \partial_2 f(x^\ast(t),x^\ast(t-1))$.

For $a\geq 0$, we denote by $({u_1}(t; a),{u_2}(t; a))$ the solution of the ODEs \eqref{eq:1.12} with the initial condition $(a, 0)$. Let us consider the amplitude derivative 
\begin{align}\label{eq:4.19}
(w_1(t),w_2(t)):= \partial_a({u_1}(t; \bar{x}), {u_2}(t; \bar{x})),
\end{align} 
which solves the linear equation \eqref{eq:4.18} with the parameter value $\beta=1$ and takes the boundary values
\begin{align}\label{eq:4.20}
\begin{split}
    (1,0)&=(w_1(0),w_2(0))\\
&=(w_1({p}),w_2({p}))+{p}'(0,\dot{x}^\ast(-1)).
\end{split}
\end{align}
Notice that $\Psi(t)$, $\dot{x}^\ast(t)$, and $w_1(t)$ are solutions of the second-order ODE
\begin{align}\label{eq:4.21}
\begin{split}
\ddot{v}_1(t)=&\left(A(t)+\frac{\dot{B}(t)}{B(t)}+A(t-1)\right)\dot{v}_1(t)\\
&+\left(\dot{A}(t)-\left(\frac{\dot{B}(t)}{B(t)}+A(t-1)\right)A(t)-\frac{1}{\beta} B(t)B(t-1)\right)v_1(t),
\end{split}
\end{align}
for parameter values $\beta=\nu$ in the case of $\Psi(t)$, and $\beta=1$ for $\dot{x}^\ast(t)$ and $w_1(t)$. Let $v^1(t)$ and $v^2(t)$ be two nonzero solutions of \eqref{eq:4.21} for parameter values $\beta_1$ and $\beta_2$, respectively. Since $B(t)\neq 0$, the only stationary solution of \eqref{eq:4.21} is $0$ and we can define the angle variables
\begin{align}
\varphi^j(t):=\arctan\left(\frac{v^j(t)}{ \dot{v}^j(t)}\right),\quad \text{for }j=1,2.\label{eq:4.22}
\end{align}
A comparison theorem \cite[Chapter 8, Theorem 1.2]{CoLe55} guarantees that for parameter values $\beta_1< \beta_2$ and initial angles $\varphi^1(0)\geq \varphi^2(0)$, the angle variables satisfy
\begin{align}\label{eq:4.23}
\varphi^1(t)> \varphi^2(t),\quad \text{for all }t>0.
\end{align}
We first prove that the initial condition $(\Psi(0),\dot{\Psi}(0))$ for the second-order ODE \eqref{eq:4.21} satisfies $\Psi(0)\neq 0$. By contradiction, we suppose that $\Psi(0)=0$ and then compare the angles
\begin{align}\label{eq:4.25}
\varphi^\mathrm{c}(t):=\arctan\left(\frac{\Psi(t)}{ \dot{\Psi}(t)}\right)\quad\text{and}\quad\varphi^\ast(t):=\arctan\left(\frac{\dot{x}^\ast(t)}{ \ddot{x}^\ast(t)}\right).
\end{align}
Since we supposed that $\Psi(0) = 0$, and we assumed $\dot x^\ast(0) = 0$, we can set $\varphi^\mathrm{c}(0)=\varphi^\ast(0)=\pi/2$. Recalling that we are in the hyperbolic setting, $\nu\neq 1$, the comparison principle \eqref{eq:4.23} yields that 
\begin{align}\label{eq:4.26}
\text{either} \quad \varphi^\mathrm{c}(t)>\varphi^\ast(t) \quad \text{or} \quad \varphi^\mathrm{c}(t)<\varphi^\ast(t), \quad \text{for all } t>0.
\end{align}

Now we recall the nodal property $z(\Psi_0)=z(\dot{x}^\ast_0)$ from Lemma \ref{Lemma4.1} and  the identity \eqref{eq:4.16}. As a result, the normalized curves
\begin{align}\label{eq:4.27}
\boldsymbol{\Psi}(t):= \frac{(\Psi(t), \dot{\Psi}(t))}{\|(\Psi(t), \dot{\Psi}(t))\|}\quad\text{and}\quad\boldsymbol{\dot{x}}(t):=\frac{(\dot{x}^\ast(t), \ddot{x}^\ast(t))}{\|(\dot{x}^\ast(t), \ddot{x}^\ast(t))\|},
\end{align}  
are both periodic with minimal period ${p}$. However, the comparison \eqref{eq:4.26} and the fact that $\boldsymbol{\Psi}(t)$ and $\boldsymbol{\dot{x}}(t)$ wind clockwise, due to the positive feedback assumption $f\in\mathfrak{X}^+$, implies that $-3\pi/2=\varphi^\ast({p})\neq \varphi^\mathrm{c}({p})= -3\pi/2$. Hence, we have reached a contradiction and $\Psi(0)\neq 0$.

Multiplying $\Psi(t)$ by a real scalar if necessary, we assume without loss of generality $(\Psi(0),\dot{\Psi}(0))=(1,\dot{\Psi}(0))$.
Now we compare the angle variable $\varphi^\mathrm{c}(t)$ to $\varphi^\zeta(t)$ given by
\begin{align}\label{eq:4.28}
\varphi^\zeta(t):= \arctan\left(\frac{\zeta(t)}{ \dot{\zeta}(t)}\right),
\end{align}
where $\zeta(t):= w_1(t)+(\dot{\Psi}(0) - \dot w_1(0))\dot{x}^\ast(t)/\ddot{x}^\ast(0)$. Here $w_1(t)$ is the first component of the amplitude derivative \eqref{eq:4.19} and therefore $\zeta(t)$ solves the second-order ODE \eqref{eq:4.21} for $\beta=1$ with the initial condition $(\zeta(0),\dot{\zeta}(0))=(\Psi(0),\dot{\Psi}(0))$. Hence, by construction, we have $\varphi^\mathrm{c}(0)=\varphi^\zeta(0)$.

Again we proceed by contradiction. Let ${p}'>0$ and suppose that $\mu_{\mathrm{c}}>1$. Then, by the inequality \eqref{eq:4.23}, it follows that
\begin{align}\label{eq:4.29}
\varphi^\zeta(t)>\varphi^\mathrm{c}(t),\quad \text{for all }t>0.
\end{align}
The normalized curve $\boldsymbol{\Psi}(t)$ winds clockwise around $(0,0)$ once in time ${p}$. By the relation \eqref{eq:4.20}, we have that $\zeta(t)$ satisfies the boundary condition
\begin{align}\label{eq:4.30}
\begin{split}
    (\zeta({p}),\dot{\zeta}({p})) &= (\zeta(0),\dot{\zeta}(0)) - {p}'(0,\ddot x^\ast(0))\\
    &=(\Psi(0),\dot{\Psi}(0)) - {p}'(0,\ddot x^\ast(0)).
\end{split}
\end{align}

Now we assemble a series of facts. First, both normalized curves $\boldsymbol{\Psi}(t)$ and 
\begin{align}\label{eq:4.31}
\boldsymbol{\zeta}(t):= \frac{(\zeta(t), \dot{\zeta}(t))}{\|(\zeta(t), \dot{\zeta}(t))\|},
\end{align} 
wind around zero clockwise. Therefore, the comparison principle \eqref{eq:4.23} implies that $\boldsymbol{\zeta}(t)$ rotates around zero faster than $\boldsymbol{\Psi}(t)$. At the same time, we use the boundary conditions \eqref{eq:4.30} to compare the endpoints $\boldsymbol{\zeta}({p})$ and $\boldsymbol{\Psi}({p})$. Second, recall that we fixed $\dot x^\ast(0)$ to be a maximum, we have that $\ddot x^\ast(0) < 0$. Since we assumed ${p}'>0$, it follows in \eqref{eq:4.30} that $-{p}'\ddot{x}^\ast(0) > 0$. In particular, $\zeta(t)$ changes signs at least twice more than $\Psi(t)$ over the time interval $[0,{p}]$. Nevertheless, $\Psi(t)$ changes signs exactly twice in that interval, and therefore $\zeta(t)$ changes signs at least four times for $t\in[0,{p}]$; see Figure \ref{fig:5}. Finally, by a comparison argument \cite[Chapter 8, Theorem 1.1]{CoLe55}, there must be a sign change of $\dot{x}^\ast(t)$ inserted between every two zeros of $\zeta(t)$. It follows that $\dot{x}^\ast(t)$ must change signs at least four times for $t\in[0,{p}]$.  By Lemma \ref{Lemma2.1} (iii), we reached a contradiction to ${p}$ being the minimal period of $\dot{x}^\ast(t)$. In particular, we have proved that if ${p}'>0$, then $\mu_{\mathrm{c}}<1$ and thereby $i(\gamma^\ast)=2n-1$.

\begin{figure}
\centering
     \begin{tabular}{Sc  Sc}
     %\hline
     {\text{$\mu_\mathrm{c} > 1,\, p' > 0$}} & {\text{$\mu_\mathrm{c} > 1,\, p' < 0$}} \\ 
      %\hline
     \begin{minipage}{0.38\textwidth}
     \cincludegraphics[width=\textwidth]{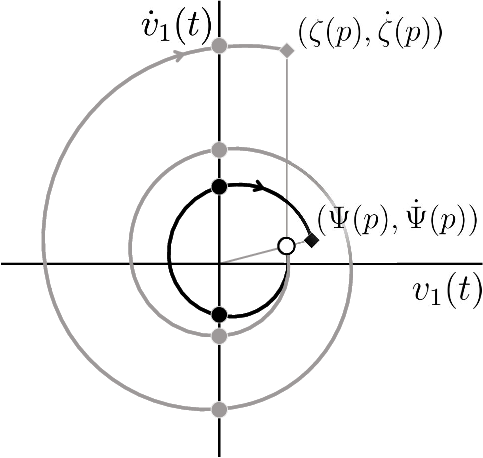}
     \end{minipage}
      & 
      \begin{minipage}{0.38\textwidth}
      \cincludegraphics[width=\textwidth]{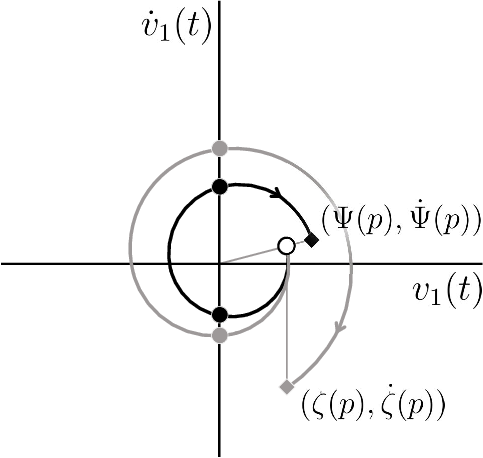}
      \end{minipage}
      \end{tabular}
\caption{Schematic plots in the $(v_1,\dot{v}_1)$-plane of the curves $(\Psi(t), \dot{\Psi}(t))$ (black) and $(\zeta(t), \dot{\zeta}(t))$ (gray) with respective angles $\varphi^\mathrm{c}(t)$ and $\varphi^\zeta(t)$, given by the formula \eqref{eq:4.25}. The curves have the same initial condition $(1, \dot{\Psi}(0))$ (white dot), the endpoint $(\Psi(p), \dot{\Psi}(p))$ (solid black diamond) corresponds to one turn around the origin as $\Psi(t)$ changes sign twice (solid black circles). (Left) The endpoint $(\zeta(p), \dot{\zeta}(p))$ (solid gray diamond) lies directly above the initial condition because $p' > 0$. Since $\mu_\mathrm{c}>1$, we have that $\varphi^\mathrm{c}(p) > \varphi^\zeta(p)$. Hence, $\zeta(t)$ must have at least four zeros (solid gray circles) for $t\in[0,p]$. (Right) The endpoint $(\zeta(p), \dot{\zeta}(p))$ (solid gray diamond) lies directly below the initial condition because $p' < 0$. Since $\mu_\mathrm{c}>1$, we have that $\varphi^\mathrm{c}(p) > \varphi^\zeta(p)$. Therefore, $(\zeta(t),\dot{\zeta}(t))$ does not complete one turn around the origin in time $p$ and $\zeta(t)$ has at most two zeros (solid gray circles) for $t\in[0,p]$.}
\label{fig:5}
\end{figure} 

To prove that  ${p}'<0$ implies $\mu_{\mathrm{c}}>1$, we follow a completely analogous argument by contradiction. In this case, however, $\zeta$ has two fewer sign changes than $\Psi$ over the time interval $[0,{p}]$. Again, this is a contradiction when we consider the sign changes of the trivial eigenfunction $\dot{x}^\ast(t)$. Therefore, ${p}'<0$ implies that $\mu_{\mathrm{c}}>1$ and the unstable dimension is $i(\gamma^\ast)=2n$.

We now go over the details of the proof in the case $f\in\mathfrak{X}^-$. The nonhyperbolic case $\mu_{\mathrm{c}}=1$ is completely analogous. If $\mu_{\mathrm{c}}\neq 1$, notice that the term $B(t)B(t-1)$ multiplying the $1/\beta$-term in the second-order ODE \eqref{eq:4.21} is positive. Therefore, the comparison principle \eqref{eq:4.23} remains unchanged. Furthermore, the criterion deciding the relative position of $(\zeta(0), \dot{\zeta}(0))$ and $(\zeta(p), \dot{\zeta}(p))$ in Figure \ref{fig:5} also remains the same. Thus, the characterization is analogous to the positive feedback case. Namely, if ${p}'<0$, then $\gamma^\ast$ has a higher unstable dimension than if ${p}'\geq 0$. Corollary \ref{Corollary4.2} (i)--(ii) yields that $i(\gamma^\ast)=2n-2$ if ${p}'>0$ and $i(\gamma^\ast)=2n-1$ otherwise.
\end{proof}

\section*{Acknowledgments}
I am deeply grateful to Bernold Fiedler and the group for Nonlinear Dynamics at the Free University of Berlin for their constant encouragement and for many insightful discussions. I would especially like to thank Jia-Yuan Dai and Babette de Wolff for carefully reading and commenting the manuscript.

\textbf{Funding.} This work has been supported by NSTC grant 113-2123-M-002-009.

\bibliographystyle{abbrv}
\bibliography{biblio} 

\begin{thebibliography}{10}

\bibitem{An88}
S.~Angenent.
\newblock The zero set of a solution of a parabolic equation.
\newblock {\em J. reine und angew. Math.}, 390, 1988.

\bibitem{BrErWa04}
P.~Brunovsk\'{y}, A.~Erdely\'{i}, and H.-O. Walther.
\newblock On a model of a currency exchange rate – local stability and
  periodic solutions.
\newblock {\em J. Dyn. Differ. Equ.}, 16, 2004.

\bibitem{ChMP78}
S.-N. Chow and J.~Mallet-Paret.
\newblock The {F}uller index and global {H}opf bifurcation.
\newblock {\em J. Differ. Equ.}, 29, 1978.

\bibitem{ChWa88}
S.-N. Chow and H.-O. Walther.
\newblock Characteristic multipliers and stability of symmetric periodic
  solutions of $\dot{x}(t)=g(x(t-1))$.
\newblock {\em Trans. Am. Math. Soc.}, 307, 1988.

\bibitem{CoLe55}
E.~A. Coddington and N.~Levinson.
\newblock {\em Theory of Ordinary Differential Equations}.
\newblock McGraw-Hill, 1955.

\bibitem{DorLa95}
P.~Dormayer and B.~Lani-Wayda.
\newblock Floquet multipliers and secondary bifurcations in functional
  differential equations: Numerical and analytical results.
\newblock {\em ZAMP}, 46, 1995.

\bibitem{FieRoWol04}
B.~Fiedler, C.~Rocha, and M.~Wolfrum.
\newblock Heteroclinic orbits between rotating waves of semilinear parabolic
  equations on the circle.
\newblock {\em J. Differ. Equ.}, 201, 2004.

\bibitem{FuRo91}
G.~Fusco and C.~Rocha.
\newblock A permutation related to the dynamics of a scalar parabolic {PDE}.
\newblock {\em J. Differ. Equ.}, 91, 1991.

\bibitem{Ha69}
J.~K. Hale.
\newblock {\em Ordinary Differential Equations}.
\newblock Krieger, 1969.

\bibitem{HaLu93}
J.~K. Hale and S.~M. Verduyn-Lunel.
\newblock {\em Introduction to Functional Differential Equations}.
\newblock Springer, 1993.

\bibitem{Hu48}
G.~E. Hutchinson.
\newblock Circular causal systems in ecology.
\newblock {\em Annals of the New York Academy of Sciences}, 50, 1948.

\bibitem{IvLaWa06}
A.~Ivanov and B.~Lani-Wayda.
\newblock Stability and instability criteria for {K}aplan-{Y}orke solutions.
\newblock {\em ZAMP}, 57, 2006.

\bibitem{KaYo74}
J.~Kaplan and J.~Yorke.
\newblock Ordinary differential equations which yield periodic solutions of
  differential delay equations.
\newblock {\em J. Math. An. and App.}, 49, 1974.

\bibitem{KaYo75}
J.~Kaplan and J.~Yorke.
\newblock On the stability of a periodic solution of a differential delay
  equation.
\newblock {\em SIMA}, 6, 1975.

\bibitem{Kr08}
T.~Krisztin.
\newblock Global dynamics of delay differential equations.
\newblock {\em Period. Math. Hung.}, 56, 2008.

\bibitem{KrGa11}
T.~Krisztin and A.~Garab.
\newblock The period function of a delay differential equation and an
  application.
\newblock {\em Period. Math. Hung.}, 63, 2011.

\bibitem{KrWa01}
T.~Krisztin and H.-O. Walther.
\newblock Unique periodic orbits for delayed positive feedback and the global
  attractor.
\newblock {\em J. Dyn. Differ. Equ.}, 13, 2001.

\bibitem{KrWaWu99}
T.~Krisztin, H.-O. Walther, and J.~Wu.
\newblock The structure of an attracting set defined by delayed and monotone
  positive feedback.
\newblock {\em CWI Quarterly}, 12, 1999.

\bibitem{LoHerKl21}
S.~Loos, S.~Hermann, and S.~Klapp.
\newblock Medium entropy reduction and instability in stochastic systems with
  distributed delay.
\newblock {\em Entropy}, 23, 2021.

\bibitem{LN23}
A.~L\'{o}pez-Nieto.
\newblock {\em Enharmonic motion: Towards the global dynamics of negative
  delayed feedback}.
\newblock PhD thesis, Freie Universität Berlin, 2023.

\bibitem{MaGla77}
M.~C. Mackey and L.~Glass.
\newblock Oscillation and chaos in physiological control systems.
\newblock {\em Science}, 197, 1977.

\bibitem{MP88}
J.~Mallet-Paret.
\newblock Morse decompositions for delay-differential equations.
\newblock {\em J. Differ. Equ.}, 72, 1988.

\bibitem{MPNu13}
J.~Mallet-Paret and R.~D. Nussbaum.
\newblock Tensor products, positive linear operators, and delay-differential
  equations.
\newblock {\em J. Dyn. Differ. Equ.}, 25, 2013.

\bibitem{MPSe962}
J.~Mallet-Paret and G.~R. Sell.
\newblock The {Poincar\'{e}-Bendixson} theorem for monotone cyclic feedback
  systems with delay.
\newblock {\em J. Differ. Equ.}, 125, 1996.

\bibitem{MPSe961}
J.~Mallet-Paret and G.~R. Sell.
\newblock Systems of differential delay equations: {Floquet} multipliers and
  discrete {L}yapunov functions.
\newblock {\em J. Differ. Equ.}, 125, 1996.

\bibitem{Ma82}
H.~Matano.
\newblock Nonincrease of the lap-number of a solution for a one-dimensional
  semilinear parabolic equation.
\newblock {\em J. Fac. Sci. Univ. Tokyo}, 1982.

\bibitem{Nuss79}
R.~Nussbaum.
\newblock Uniqueness and nonuniqueness for periodic solutions of
  $\dot{x}(t)=-g(x(t-1))$.
\newblock {\em J. Differ. Equ.}, 34, 1979.

\bibitem{Nuss73}
R.~D. Nussbaum.
\newblock Periodic solutions of some nonlinear, autonomous functional
  differential equations. {II}.
\newblock {\em J. Differ. Equ.}, 14, 1973.

\bibitem{SuaSch88}
M.~Suarez and P.~Schopf.
\newblock A delayed action oscillator for enso.
\newblock {\em J. Atmos. Sci.}, 45, 1988.

\bibitem{Van82}
A.~Vanderbauwhede.
\newblock {\em Local bifurcation and symmetry}.
\newblock Pitman, 1982.

\bibitem{LaWa76}
M.~Wazewska-Czyzewska and A.~Lasota.
\newblock Mathematical problems of the dynamics of red blood cells system.
\newblock {\em Ann. Pol. Math. Soc.}, 4, 1976.

\bibitem{YoKo20}
K.~Yoshioka-Kobayashi, M.~Matsumiya, Y.~Niino, A.~Isomura, H.~Kori,
  A.~Miyawaki, and R.~Kageyama.
\newblock Coupling delay controls synchronized oscillation in the segmentation
  clock.
\newblock {\em Nature}, 580, 2020.

\end{thebibliography}

\end{document}